\documentclass[12pt]{article}

\usepackage{verbatim}

\usepackage{amsfonts}

\usepackage{amsthm}

\usepackage{amssymb}

\usepackage{amsmath}

\usepackage{mathabx}

\usepackage{enumerate}

\usepackage[all]{xy}

\usepackage{graphicx}

\usepackage{hyperref}

\newcommand{\N}{\mathbb{N}}

\newcommand{\Z}{\mathbb{Z}}

\newcommand{\R}{\mathbb{R}}

\newcommand{\C}{\mathbb{C}}

\newcommand{\Hawaii}{Hawai\kern.05em`\kern.05em\relax i}

\newcommand{\ra}{~~\Rightarrow~~}

\numberwithin{equation}{section}

\theoremstyle{plain}
\newtheorem{theorem}[equation]{Theorem}
\newtheorem{lemma}[equation]{Lemma}
\newtheorem{corollary}[equation]{Corollary}
\newtheorem{proposition}[equation]{Proposition}

\newtheorem{definition-theorem}[equation]{Definition / Theorem}

\newtheorem*{conjecture*}{Conjecture}
\newtheorem*{theorem*}{Theorem}

\theoremstyle{definition}
\newtheorem{definition}[equation]{Definition}
\newtheorem{example}[equation]{Example}

\theoremstyle{remark}
\newtheorem{remark}[equation]{Remark}

\newtheorem*{example*}{Example}  
\newtheorem*{remark*}{Remark}

\title{Dynamic Asymptotic Dimension: relation to dynamics, topology, coarse geometry, and $C^*$-algebras}

\author{Erik Guentner, Rufus Willett, and Guoliang Yu}

\begin{document}

\maketitle

\begin{abstract}
We introduce dynamic asymptotic dimension, a notion of dimension for actions of discrete groups on locally compact spaces, and more generally for locally compact \'{e}tale groupoids.  We study our notion for minimal actions of the integer group, its relation with conditions used by Bartels, L\"{u}ck, and Reich in the context of controlled topology, and its connections with Gromov's theory of asymptotic dimension.  We also show that dynamic asymptotic dimension gives bounds on the nuclear dimension of Winter and Zacharias for $C^*$-algebras associated to dynamical systems.  Dynamic asymptotic dimension also has implications for $K$-theory and manifold topology: these will be drawn out in subsequent work.
\end{abstract}

\thanks{The first author was partially supported by a grant from the
  Simons Foundation (\#245398).
  The second author was partially supported by NSF grant DMS-1401126.  The third author was partially supported by NSF grant DMS-1362772.}

\tableofcontents

\section{Introduction}

The main aim of this paper is to introduce \emph{dynamic asymptotic dimension}, a property of topological dynamical systems.  Precisely, we are interested in actions of discrete groups on locally compact spaces: throughout the paper, we say `$\Gamma\lefttorightarrow X$ is an action' as shorthand for saying that $\Gamma$ is a discrete group acting by homeomorphisms on a locally compact Hausdorff topological space $X$.

To give an idea of our main definition, we will state it here in a specialized form.  First, we need a preliminary on `broken' orbit equivalence relations.

\begin{definition}
Let $\Gamma\lefttorightarrow X$ be an action.  For a subset $E$ of $\Gamma$ and an open subset $U$ of $X$, let $\sim_{U,E}$ be the equivalence relation on $U$ generated by $E$: precisely, for $x,y\in U$, $x\sim_{U,E} y$ if there is a finite sequence 
$$
x=x_0,x_1,...,x_n=y
$$
of points in $U$ such that for each $i\in \{1,...,n\}$ there exists $g\in E\cup E^{-1}\cup \{e\}$ such that $gx_{i-1}=x_i$.  
\end{definition}

Note that if $U=X$ then $\sim_{U,E}$ is just the equivalence relation of being in the same orbit for the subgroup $\langle E\rangle$ of $\Gamma$ generated by $E$.  However, if $U$ is a proper subset of $X$ then $\sim_{U,E}$ equivalence classes will generally be smaller than the intersection of $U$ with the $\langle E\rangle$-orbits.

Here is our main definition.

\begin{definition}\label{firstdad}
Let $\Gamma\lefttorightarrow X$ be an action, where we assume for simplicity that $X$ is compact, and the action is free\footnote{This means that if $gx=x$ for some $g\in \Gamma$ and $x\in X$, then $g=e$ is the identity element of $\Gamma$.}.
The \emph{dynamic asymptotic dimension} of $\Gamma\lefttorightarrow X$ is the smallest $d\in \N$ with the following property: for each finite subset $E$ of $\Gamma$, there is an open cover $\{U_0,...,U_d\}$ of $X$ such that for each $i\in \{0,...,d\}$, the equivalence relation $\sim_{U_i,E}$ on $U_i$ has uniformly finite equivalence classes.
\end{definition}

Heuristically, we think that $\Gamma\lefttorightarrow X$ has dynamic asymptotic dimension at most $d$, if for any (large) finite subset $E$ of $\Gamma$, the action can be `broken' into at most $d+1$ parts, and on each part the `action generated by $E$' has only `finite complexity'.  We generalize this definition to non-free actions on non-compact spaces in Definition \ref{daddef} below, and then to locally compact, Hausdorff, \'{e}tale groupoids in Definition \ref{gpddad}.

Our main motivation for introducing this property is its implications for $K$-theory of associated algebras (and more general categories) and thus for manifold topology.  These implications come via the use of controlled cutting-and-pasting, or Mayer-Vietoris, techniques pioneered by the third author \cite{Yu:1998wj} in the setting of asymptotic dimension, and developed in a more general context by the first and third authors in collaboration with Tessera \cite{Guentner:2009tg}.  Other important motivations come from work of Farrell and Jones \cite{Farrell:1986ve}, of Bartels, L\"{u}ck and Reich \cite{Bartels:2008uq}, and of Bartels and L\"{u}ck \cite{Bartels:2012fu} in controlled topology.  We will explore these aspects in other work \cite{Guentner:2014bh}.

We believe however, that dynamic asymptotic dimension will admit many interesting examples, and be useful in other contexts: it is the purpose of this paper to explore some of these other aspects.  Specifically, we will develop some of the main motivating examples and discuss some consequences for the structure theory of $C^*$-algebras.

\subsection*{Examples}

Our main theorems on examples are as follows.  

\begin{theorem}\label{exthe}
\begin{enumerate}[(i)]
\item Let $\Z\lefttorightarrow X$ be a free, minimal $\Z$ action on a compact space.  Then the dynamic asymptotic dimension of $\Z\lefttorightarrow X$ is one.
\item Let $\Gamma\lefttorightarrow X$ be an action satisfying a `Bartels-L\"{u}ck-Reich type condition of dimnesion $d$' for the family of finite subgroups.  Then the dynamic asymptotic dimension of $\Gamma\lefttorightarrow X$ is at most $d$.  
\item The canonical action $\Gamma\lefttorightarrow \beta\Gamma$ of any (countable) discrete group on its Stone-\v{C}ech compactification is equal to the asymptotic dimension of $\Gamma$ in the sense of Gromov.  More generally, if $X$ is a bounded geometry coarse space and $G(X)$ the associated coarse groupoid, then the dynamic asymptotic dimension of $G(X)$ equals the asymptotic dimension of $X$.
\item Let $\Gamma$ be a countable group with asymptotic dimension $d$.  Then $\Gamma$ admits a free minimal action on the Cantor set with dynamic asymptotic dimension at most $d$.
\end{enumerate}
\end{theorem}

See Theorems \ref{minz}, \ref{blrlp}, \ref{fadthe} and Corollary \ref{mincor} below for parts (i), (ii), (iii) and (iv) respectively (and more explanation of the terminology involved in each part).  

Part (i) says that the actions that are perhaps most interesting from the point of view of classical topological dynamics fall under the purview of dynamic asymptotic dimension in a natural and simple way.  On the other hand, part (iv) implies that many interesting classes of groups---for example, word hyperbolic groups \cite{Roe:2005rt}, CAT(0) cubical groups \cite{Wright:2012aa}, lattices in Lie groups, and many solvable groups \cite{Bell:2008fk}---admit at least some actions with finite dynamic asymptotic dimension.  Parts (ii) and (iii) were our principal motivations, as they give the connections to controlled $K$-theory and controlled topology which underlie our work in those directions.

\subsection*{Implications}

The main implications we explore in this paper are to the structure theory of $C^*$-algebras: in particular to nuclear dimension in the sense of Winter and Zacharias \cite{Winter:2010eb}, a property that has been very important in Elliott's classification program \cite{Elliott:1995dq} and elsewhere.  The main result we have here is Theorem \ref{fnd} below: this says that under a minor technical hypothesis, the nuclear dimension of the reduced $C^*$-algebra of a free \'{e}tale groupoid can be bounded in terms of the dynamic asymptotic dimension of the groupoid, and the covering dimension of the groupoid's unit space.  Rather than repeat Theorem \ref{fnd} here, we just give some corollaries.  See Section \ref{fndsec} for more details.

\begin{theorem}\label{impthe}
\begin{enumerate}[(i)]
\item Let $\Z\lefttorightarrow X$ be a free, minimal action on a second countable compact space of covering dimension $N$.  Then the nuclear dimension of the crossed product $C(X)\rtimes \Z$ is at most $2N+1$.
\item Let $X$ be a bounded geometry coarse space.  Then the nuclear dimension of the uniform Roe algebra $C^*_u(X)$ is at most the asymptotic dimension of $X$.
\item Any countable group $\Gamma$ admits a free, minimal action on the Cantor set $X$ such that the associated reduced crossed product $C(X)\rtimes_r\Gamma$ has nuclear dimension at most the asymptotic dimension of $\Gamma$.
\end{enumerate}
\end{theorem}

Parts (i) and (ii) are originally due to Toms and Winter \cite[Section 3]{Toms:2013vn} and Winter and Zacharias \cite[Section 8]{Winter:2010eb} respectively; moreover, our proofs are in some sense close to the original ones (if heavily disguised).  Nonetheless, we think there is some interest in explicitly bringing these results under one dynamical framework.  The final result seems to be new, and says that `many' groups admit simple crossed products with finite nuclear dimension.  Thanks to spectacular recent advances in $C^*$-algebra theory \cite{Elliott:2015fb,Tikuisis:2015kx}, it seems that one now knows that this implies these crossed products fall under the purview of the Elliott program.

\subsection*{Outline of the paper}

In Section \ref{dadsec} we introduce our main definition in the case of group actions, and mention some basic consequences.  In Section \ref{minzsec} we study minimal $\Z$ actions on compact spaces using ideas of Putnam \cite{Putnam:1989hi}.  In Section \ref{blrsec} we discuss the connection to the work of Bartels, L\"{u}ck, and Reich \cite{Bartels:2008uq,Bartels:2012fk} in controlled topology.  In Section \ref{gpdsec}, we extend the main definition to \'{e}tale groupoids with applications to coarse geometry in the next section in mind.  In Section \ref{asdimsec} we explore the connections to coarse geometry, and use this to construct the examples in part (iv) of Theorem \ref{exthe}, as well as to clarify the relationship of the Bartels-L\"{u}ck-Reich conditions to asymptotic dimension.  Section \ref{pousec} is devoted to a technical construction of almost invariant partitions of unity.  Finally, Section \ref{fndsec} discusses the implications to nuclear dimension of groupoid $C^*$-algebras, and in particular all the parts of Theorem \ref{impthe}.

We have tried to write the paper in a `modular' way, so that sections can be read independently of each other to a large extent.  In particular, groupoids are not mentioned until Section \ref{gpdsec}, and (noncommutative) $C^*$-algebras are not mentioned before Section \ref{fndsec}.

\subsection*{Acknowledgements}

This paper has been some time in gestation, and has benefited from conversations with several people.  In particular, we would like to thank Arthur Bartels, Siegfried Echterhoff, David Kerr, Ian Putnam, Daniel Ramras, Wilhelm Winter, and Jianchao Wu for useful comments, penetrating questions, and / or patient explanations.  The first two authors would like to thank Texas A\&{}M University and the Shanghai Center for Mathematical Sciences for their hospitality during some of the work on this paper.

\section{Dynamic asymptotic dimension for group actions}\label{dadsec}

In this section, we introduce our main definition.

\begin{definition}\label{daddef}
An action $\Gamma\lefttorightarrow X$ has \emph{dynamic asymptotic dimension} $d$ if $d$ is the smallest natural number with the following property: for any compact subset $K$ of $X$ and finite subset $E$ of $\Gamma$ there are open subsets $U_0,...,U_d$ of $X$ that cover $K$ such that for each $i\in \{0,...,d\}$, the set 
\begin{equation}\label{finsets}
\left\{\begin{array}{l|l}  & \text{there exist } x\in U_i\text{ and }g_n,...,g_1\in E\text{ such that } \\ g\in \Gamma &g=g_n\cdots g_2g_1 \text{ and  for all } k\in \{1,...,n\}, \\ &  g_k\cdots g_1x\in U_i\end{array}\right\}
\end{equation}
is finite. 

The action has \emph{finite dynamic asymptotic dimension} if it has dynamic asymptotic dimension $d$ for some $d$, and \emph{infinite dynamic asymptotic dimension} otherwise.
\end{definition}

\begin{remark}\label{dadrem}
\begin{enumerate}[(i)]
\item \label{dad0} An action $\Gamma\lefttorightarrow X$ on a compact space $X$ has dynamic asymptotic dimension $0$ if and only if $\Gamma$ is locally finite, i.e.\ any finite subset of $\Gamma$ generates a finite subgroup. 
\item \label{dad02} Recall that an action $\Gamma\lefttorightarrow X$ is \emph{proper} if for any compact subset $K$ of $X$, the set $\{g\in \Gamma~|~gK\cap K\neq \varnothing\}$ is finite.  It is \emph{locally proper} if every finitely generated subgroup of $\Gamma$ acts properly on $X$.  Point \eqref{dad0} says that an action on a compact space has dynamic asymptotic dimension zero if and only it is locally proper.

In general, it is still true that local properness implies dynamic asymptotic dimension zero: indeed, let a finite $E\subseteq \Gamma$ and compact $K\subseteq X$ be given.  Let $U_0$ be any relatively compact open set that contains $K$ (such a $U_0$ exists as $X$ is locally compact).  Then if $\langle E\rangle$ is the subgroup of $\Gamma$ generated by $E$, the set in line \eqref{finsets} above is contained in 
$$
\{g\in \langle E\rangle ~|~g^{-1}U_0\cap U_0\neq \varnothing \}
$$
and is thus finite by local properness.  The converse is false, however: for example, the action of $\Z$ on $\R^2$ considered in \cite[Chapter 5, Example 8.3]{Massey:1977sf} has dynamic asymptotic dimension zero but is not (locally) proper.  It is perhaps most natural to characterize dynamic asymptotic dimension zero actions in terms of groupoids: see Example \ref{lfgpd} below.
\item If an action has finite dynamic asymptotic dimension, then all point stablizers must be locally finite.  In particular, an action of a torsion free group with finite dynamic asymptotic dimension is free.
\item If $\Gamma\lefttorightarrow X$ is a free action with $X$ compact, then Definition \ref{daddef} above is equivalent to Definition \ref{firstdad} from the introduction.  

To see this, assume first that $\Gamma\lefttorightarrow X$ satisfies the condition in Definition \ref{daddef}, and let a finite subset $E$ of $\Gamma$ be given.  Take $K=X$, and a cover $U_0,...,U_d$ as in Definition \ref{daddef} for this $K$ and $E$.  Then if $F_i$ is the finite set in line \eqref{finsets}, all $\sim_{U_i,E}$ equivalence classes are contained in $F_i\cdot x$ for some $x$, and thus have uniformly bounded finite cardinalities.

Conversely say a finite subset $E$ of $\Gamma$ is given and $\{U_0,...,U_d\}$ is a cover of $X$ with the properties in Definition \ref{firstdad}, so in particular all $\sim_{U_i,E}$ equivalence classes have cardinality at most some integer $N$.  Then each of the sets $F_i$ in line \eqref{finsets} is contained in $(E\cup E^{-1}\cup\{e\})^N$, and thus has cardinality at most $(2|E|+1)^N$.
\item See Corollary \ref{dadgppou} below for a stronger-looking equivalent definition of dynamic asymptotic dimension in terms of almost-invariant partitions of unity.
\end{enumerate}
\end{remark}

\section{Example - minimal $\Z$-actions}\label{minzsec}

In this section we give our first non-trivial examples: minimal actions of $\Z$ on compact spaces (recall that an action $\Gamma\lefttorightarrow X$ is \emph{minimal} if all $\Gamma$-orbits are dense).   

The proof is inspired by ideas of Ian Putnam that were originally used to build interesting AF-algebras associated to minimal actions of $\Z$ on the Cantor set \cite{Putnam:1989hi} (and subsequently much developed by Putnam and others).  We would like to thank Prof. Putnam for suggesting that there might be a connection between his work and our notion of dynamic asymptotic dimension.

\begin{theorem}\label{minz}
Let $\Z\lefttorightarrow X$ be a minimal $\Z$ action on an infinite compact space $X$.  Then the dynamic asymptotic dimension of the action is one.
\end{theorem} 

\begin{proof}
As $\Z$ is not locally finite, Remark \ref{dadrem} part \eqref{dad0} implies that the dynamic asymptotic dimension of $\Z\lefttorightarrow X$ is not zero, so it suffices to bound the dynamic asymptotic dimension above by one.

Let a finite subset $E$ of $\Z$ be given; we may as well assume that $E$ is an `interval' $[-N,N]\cap \Z$ for some $N\in \N$.  As the action is minimal and $X$ is infinite, the action is free.  It follows that we can find a non-empty open subset $U$ of $X$ such that 
$$
n\cdot U\cap U=\varnothing
$$
for all $n\in [-5N,5N]\setminus \{0\}$.  Let $V$ be any non-empty open subset of $X$ such that $\overline{V}\subseteq U$.  Define
$$
U_0:=\bigcup_{n=-N}^N n\cdot U,~~~U_1:=X~\setminus~ \bigcup_{n=-N}^N n\cdot \overline{V}.
$$
Clearly $\{U_0,U_1\}$ is an open cover of $X$.  To finish the proof, it suffices to show that the sets 
\begin{equation}\label{finsetsz}
\left\{\begin{array}{l|l}  & \text{there exist } x\in U_i\text{ and }n_1,...,n_m\in E\text{ such that } \\ n\in \Z &n=n_1+\cdots +n_m  \text{ and  for all } k\in \{1,...,m\}, \\ &  (n_k+\cdots+n_1)x\in U_i\end{array}\right\}
\end{equation}
are finite for $i=0,1$.\\

First look at $U_0$.  We claim that in this case the set in line \eqref{finsetsz} is contained in $[-3N,3N]$.  Indeed, assume for contradiction that there exist $n_1,...,n_m\in [-N,N]$ with $|n_1+\cdots + n_m|>3N$ and $x\in U_0$ with $(n_k+\cdots +n_1)x\in U_0$ for all $k\in \{1,...,m\}$.  Write $x=n_0x_0$ for some $x_0\in U$ and $n_0\in [-N,N]$ (such exist by definition of $U_0$).   As each $n_k$ is in $[-N,N]$ and $|n_m+\cdots +n_1|>3N$, there must exist $k$ such that $3N\geq |n_k+\cdots +n_1|>2N$, whence 
\begin{equation}\label{bound1}
4N\geq |n_k+\cdots +n_1+n_0|>N.
\end{equation}
As 
$$
(n_k+\cdots +n_1)x=(n_k+\cdots +n_1+n_0)x_0
$$
is in $U_0$, there are $n_0'\in [-N,N]$ and $x_0'\in U$ with 
$$
(n_k+\cdots +n_1+n_0)x_0=n_0'x_0',
$$
whence
$$
(n_k+\cdots +n_1+n_0-n_0')x_0=x_0'
$$
Hence as $x_0$ and $x_0'$ are in $U$ we have that
$$
(n_k+\cdots +n_1+n_0-n_0')\cdot U\cap U\neq \varnothing.
$$
However, line \eqref{bound1} implies that $5N\geq |n_k+\cdots +n_1+n_0-n_0'|>0$, so this contradicts that $U\cap n\cdot U=\varnothing$ for $n\in [-5N,5N]\setminus \{0\}$.\\

We now look at the set in line \eqref{finsetsz} for $U_1$.  We first claim that there exists $M\in \N$ such that for all $x\in U_1$ there exists $m_-\in [-M,0)$ and $m_+\in (0,M]$ with $m_-x\in V$ and $m_+x\in V$.  Indeed, for each $M\in \N$ define
$$
W_M:=\left\{\begin{array}{l|l} x\in X & \text{there are } m_-\in [-M,0),m_+\in (0,M]  \\ & \text{ such that } m_-x,m_+x\in V\end{array}\right\}.
$$
It follows from the fact that $V$ is open that each $W_M$ is open.   Moreover, minimality of the action implies that for each $x\in X$, the `half-orbits' $\N\cdot x$ and $(-\N)\cdot x$ are dense: if not, the limit points of one of these sets would be a closed $\Z$-invariant subset.  Hence the fact that $V$ is open implies that each $x\in X$ is in $W_M$ for some $M$, and so $\{W_M~|~M\in \N\}$ is an open cover of $X$.   Compactness of $X$ and the fact that $W_{M_1}\subseteq W_{M_2}$ for $M_1\leq M_2$ implies that $X$ is contained in $W_M$ for some $M$, and this implies the claim.

To complete the proof, we will now show that for $U_1$, the set in line \eqref{finsetsz} is contained in $[-M-N,M+N]$.  Assume for contradiction this fails, so there exist $n_1,...,n_m\in E$ with $|n_1+\cdots+n_m|>M+N$ and $x\in U_1$ such that $(n_1+\cdots +n_k)x$ is in $U_1$ for all $k\in \{1,...,m\}$.  Assume for simplicity that $n_1+\cdots +n_m>M+N$; the case $n_1+\cdots n_m<-M-N$ can be handled similarly.  Let $m_+\in (0,M]$ be such that $m_+x\in V$.  Then there exists $k\in \{1,...,m\}$ with $(n_1+\cdots +n_k)-m_+\in [-N,N]$.  In particular, 
$$
(n_1+\cdots +n_k)x=(m_++n)x
$$
for some $n\in [-N,N]$.  However as $m_+x$ is in $V$, $(m_++n)x$ is in $X\setminus U_1$ by definition of $U_1$, which is a contradiction.

\end{proof}

\section{Example - Bartels-L\"{u}ck-Reich conditions}\label{blrsec}

Our main goal in this section is to study some properties of group actions that were important in the work of Bartels, L\"{u}ck, and Reich \cite{Bartels:2008uq,Bartels:2012fk} on the Farrell-Jones conjecture and show that that they imply that the action has finite dynamic asymptotic dimension.  

We will need to establish come conventions on simplicial complexes.  Let $V$ be a set, thought of a discrete topological space.  The \emph{space of probability measures} on $V$ is
$$
P(V):=\{\mu\in l^1(V)~|~\mu(v)\geq 0 \text{ for all }v\in V \text{ and } \|\mu\|_1=1\},
$$ 
equipped with the metric 
$$
d\Big(\sum_{v\in V}t_vv~,~\sum_{v\in V}s_vv\Big):=\sum_{v\in V}|t_v-s_v|.
$$
coming from the $l^1$-norm.   Write $P_n(V)$ for the subset of $P(V)$ consisting of measures supported on at most $n+1$ points ad define
$$
P_f(V):=\bigcup_{n=0}^\infty P_n(V)
$$ 
to be the subspace of finitely supported probability measures. Elements of $P(V)$ will usually be written as formal sums
$$
\mu=\sum_{v\in V}t_vv,
$$
where $\mu(v)=t_v$ is in $[0,1]$ and $\sum_vt_v=1$.  We identify $V$ with $P_0(V)$ in the obvious way.

A \emph{simplicial complex} is a subspace $C$ of some $P_f(V)$ as above such that $V\subseteq C$, and with the property that whenever 
$$
\mu=\sum_{v\in V}t_vv
$$ 
is an element of $C$ and $S=\{v\in V~|~t_v\neq 0\}$ is the support of $\mu$, then 
$$
\Big\{\sum_{v\in V}t_vv\in P_f(V)~\Big|~t_v=0 \text{ for } v\not\in S\Big\}
$$
is contained in $C$.  The \emph{vertex set} of $C$ is $V$.  A simplicial complex $C$ is equipped with the restriction of the $l^1$-metric defined above.   If $C$ is a simplicial complex and $n\in \N$, an \emph{$n$-simplex} in $C$ is a subset that is equal to the convex hull of some set of $n+1$ vertices. The \emph{$n$-skeleton of $C$} is defined to be $C_n:=C\cap P_{n}(V)$ (so in particular $C_0=V$ and $C_{-1}=\varnothing$).  The \emph{dimension} of $C$ is the smallest $d$ such that $C=C_d$ (or infinity if no such $d$ exists).

If $\Gamma$ is a discrete group and $C$ a simplicial complex with vertex set $C$, then a \emph{simplicial action} of $\Gamma$ on $C$ is an action that is induced from some action of $\Gamma$ on $V$ via the formula
$$
g\Big(\sum_{v\in V}t_vv\Big)=\sum_{v\in V}t_v(gv).
$$
Note that a simplicial action is isometric.  All actions on simplicial complexes will be assumed simplicial.  

We will need the following technical lemma at a couple of points below: roughly, it says that simplicial complexes admit covers with rather rigid combinatorial properties.

\begin{lemma}\label{nice cover}
Let $C$ be a simplicial complex of dimension at most $d$ and equipped with a simplicial action by $\Gamma$.  For a subset $A$ of $C$ and $\delta>0$, write 
$$
N_\delta(A):=\{x\in C~|~d(x,A)<\delta\}
$$
for the $\delta$-neighbourhood of $A$.  For each $i\in \{0,...,d\}$ define
$$
V_i:=N_{\frac{1}{3}10^{-i}}(C_i)~\setminus ~\overline{N_{\frac{5}{2}10^{-i}}(C_{i-1})}.
$$
Then the collection $\{V_0,...,V_d\}$ is an open cover of $C$ by $\Gamma$-invariant subsets.  

Moreover, for each $i\in \{0,...d\}$ and $i$-simplex $\Delta$, define 
$$
V_{i~\Delta}:=N_{\frac{1}{3}10^{-i}}(\Delta)~\setminus ~\overline{N_{\frac{5}{2}10^{-i}}(C_{i-1})}.
$$
Then 
$$
V_i=\bigcup_{\Delta \text{ an $i$-simplex}}V_{i~\Delta},
$$
the $\Gamma$ action permutes the distinct $V_{i~\Delta}$, and for distinct $i$-simplices $\Delta$ and $\Delta'$, 
\begin{equation}\label{separation}
d(V_{i~\Delta},V_{i~\Delta'})\geq \frac{1}{3}10^{-i}\geq \frac{1}{3}10^{-d}.
\end{equation}
\end{lemma}

\begin{proof}
Clearly each $V_i$ is open, and each is $\Gamma$-invariant as $\Gamma$ acts simplicially (whence it preserves each $i$-skeleton and the metric).  To show that $\{V_0,...,V_d\}$ covers $C$, it suffices to show that $C_i\subseteq V_0\cup\cdots \cup V_i$ for each $i$, which is clear by induction.  The decomposition $V_i=\bigcup_{\Delta \text{ an $i$-simplex}}V_{i~\Delta}$ is also clear as $C_i=\bigcup_{\Delta \text{ an $i$-simplex}}\Delta$.  Moreover, as the action is simplicial $\Gamma$ permutes the sets $V_{i~\Delta}$.

Finally, say for contradiction that $\Delta$, $\Delta'$ are distinct $i$-simplices such that $d(V_{i~\Delta},V_{i~\Delta'})<\frac{1}{3}10^{-i}$.  Then there exists $\mu\in V_{i~\Delta}$ such that $d(\mu,\Delta')<\frac{2}{3}10^{-i}$.  Note that
$$
\frac{1}{3}10^{-i}>d(\mu,\Delta)\geq \sum_{v\not\in \Delta}t_v,\quad\text{and}\quad \frac{2}{3}10^{-i}>d(\mu,\Delta')\geq \sum_{v\not\in \Delta'}t_v
$$
whence
$$
T:=\sum_{v\not\in \Delta\cap \Delta'}t_v\leq\sum_{v\not\in \Delta}t_v+ \sum_{v\not\in \Delta'}t_v <10^{-i}.
$$
Hence in particular $T<1$ and so we may define
$$
\nu:=\sum_{v\in \Delta\cap \Delta'}\frac{t_v}{1-T}v,
$$
which is in $C_{i-1}$ as $\Delta$ and $\Delta'$ are distinct.  Note that 
$$
d(\mu,C_{i-1})\leq d(\mu,\nu)=\sum_{v\in \Delta\cap \Delta'}t_v\Big(\frac{1}{1-T}-1\Big)+\sum_{v\not\in \Delta\cap \Delta'}t_v=2T<2\cdot 10^{-i}.
$$
However, as $\mu$ is in $V_{i~\Delta}$ we have $\mu\not\in \overline{N_{\frac{5}{2}10^{-i}(C_{i-1})}}$, and putting this together with the line above gives
$$
\frac{5}{2}10^{-i}<d(\mu,C_{i-1})<2\cdot 10^{-i},
$$
which is the desired contradiction.
\end{proof}

A \emph{family} of subgroups of a discrete group $\Gamma$ is a collection of subgroups satisfying the following conditions:
\begin{itemize}
\item $\mathcal{F}$ is closed under conjugation;
\item $\mathcal{F}$ is closed under taking subgroups;
\item $\mathcal{F}$ is closed under taking finite index supergroups.
\end{itemize}
The examples that are important for applications are the family of finite subgroups, and that of virtually cyclic subgroups.   

\begin{definition}\label{g-f complex}
Let $C$ be a simplicial complex equipped with a simplicial action of a discrete group $\Gamma$.  Let $\mathcal{F}$ be a family of subgroups of $\Gamma$ satisfying the conditions above.  Then $C$ is called a \emph{$(\Gamma,\mathcal{F})$-complex} if the stabilizer of every vertex in $C$ is an element of $\mathcal{F}$.
\end{definition}

We will only use the following definition in the case that $Y$ is a simplicial complex equipped with a simplicial action of $\Gamma$.

\begin{definition}\label{aemaps}
Let $X$ be a topological space and $Y$ a metric space, and assume that $\Gamma$ is a discrete group acting on $X$ by homeomorphisms, and on $Y$ by isometries.  Let $E$ be a subset of $\Gamma$ and $\epsilon>0$.  A map $f:X\to Y$ is \emph{$(E,\epsilon)$-equivariant} if 
$$
\sup_{x\in X}d_Y(f(g x),g f(x))<\epsilon.
$$
for all $g\in E$.
\end{definition}

As the last of our preliminaries before getting to the conditions of Bartels-L\"{u}ck-Reich, we have another technical lemma.

\begin{lemma}\label{perturb}
Let $f:X\to C$ be a continuous $(E,\epsilon)$-equivariant map as in Definition \ref{aemaps}, where $X$ is compact and $E$ is finite.  Then there exists a finite subset $S$ of the vertex set of $C$ and an $(E,\epsilon)$-equivariant map $f':X\to C$ such that $f'(X)\subseteq P(S)\cap C$.
\end{lemma}

\begin{proof}
Note that for any fixed $\delta>0$, the collection 
$$
\{N_\delta(P(S)\cap C)~|~S\text{ a finite set of vertices}\}
$$
is an open cover of $C$.  It follows that if $K\subseteq C$ is compact, then for any $\delta>0$ there exists a finite subset $S$ of $V$ such that $K\subseteq N_\delta(P(S)\cap C)$.  Define
$$
\delta:=\min\Big\{1~,~\frac{1}{2}\min_{g\in E}\big\{\epsilon-\sup_{x\in X}d_C(f(g x),g f(x))\big\}\Big\}>0
$$
and let $S$ be a finite subset of $V$ such that $f(X)\subseteq N_{\delta/2}(P(S)\cap C)$.  
We may write $f(x)=\sum_{v\in V}t_v(x)v$ where each $t_v:X\to[0,1]$ is a continuous function.  Note that for any $\mu=\sum_{v\in S} s_v v\in P(S)$,
$$
d(\mu,f(x))=\sum_{v\in S}|s_v-t_v(x)|+\sum_{v\not\in S}t_v(x)\geq \sum_{v\not\in S}t_v(x).
$$
Taking the infimum over all such $\mu$ and using that $f(X)\subseteq N_{\delta/2}(P(S)\cap C)$ gives 
$$
\sum_{v\not\in S}t_v(x)<\delta/2.
$$
Hence in particular the formula $T(x):=\sum_{v\in S}t_v(x)$ defines a continuous function $T:X\to (1-\delta/2,1]$.  Define 
$$
f'(x):=\sum_{v\in S}\frac{t_v(x)}{T(x)}v,
$$
so $f':X\to C\cap P(S)$ is a continuous function.  Then for any $x\in X$,
$$
d_C(f(x),f'(x))=\sum_{v\in S}t_v(x)\Big(\frac{1}{T(x)}-1\Big)+\sum_{v\not\in S}t_v(x)=2(1-T(x))<\delta.
$$
It follows that $f':X\to C$ has the desired properties.
\end{proof}

We now come to the conditions that are our main object of study in this section.  The first condition in the following proposition is essentially taken from Bartels' survey paper  \cite[Theorem A, page 9]{Bartels:2012fk}, and the second from the paper of Bartels-L\"{u}ck-Reich on equivariant covers of hyperbolic groups \cite[Theorem 1.2]{Bartels:2008uq}.  The result is very closely connected to \cite[Lemma 4.4]{Szabo:2014aa}, and is to some extent already implicit in \cite{Bartels:2008uq}; nonetheless, we do not think a complete proof exists anywhere in the literature, so give one here.

\begin{proposition}\label{blr}
Say $\Gamma\lefttorightarrow X$ is an action with $X$ compact, and $\mathcal{F}$ is a family of subgroups of $\Gamma$.  The following are equivalent.
\begin{enumerate}[(i)]
\item There exists $d$ such that for all finite $E\subseteq\Gamma$ and all $\epsilon>0$ there exists a $(\Gamma,\mathcal{F})$-complex $C$ of dimension at most $d$ and and a continuous $(E,\epsilon)$-equivariant map
$$
f:X\to C.
$$
\item There exists $d$ such that for each finite subset $E$ of $\Gamma$ there exists a $\Gamma$ equivariant open cover $\mathcal{U}$ of $X\times\Gamma$  such that the following hold.
\begin{enumerate}[(A)]
\item For every $U\in\mathcal{U}$ and $g\in \Gamma$, $g U=U$, or $gU\cap U=\varnothing$.
\item For every $U\in\mathcal{U}$, $\{g\in \Gamma~|~gU=U\}\in \mathcal{F}$.
\item The multiplicity of $\mathcal{U}$ is at most $d+1$.
\item $\mathcal{U}/\Gamma$ is finite.
\item For every $g\in \Gamma$ and $x\in X$ there exists $U\in\mathcal{U}$ such that $\{x\}\times gE\subseteq U$. 
\end{enumerate}
\end{enumerate}
\end{proposition}

\begin{proof}
Assume condition (i).  Let $\epsilon=\frac{1}{6}10^{-d}$.  Fix a finite subset $E$ of $\Gamma$, and let $C$, $f$ be as in condition (i) for this choice of $E$ and $\epsilon$.  Using Lemma \ref{perturb}, we may assume that the image of $f$ is contained in $P(S)\cap C$ for some finite subset $S$ of the vertex set of $C$.  Define
$$
\phi:X\times\Gamma\to C\cap P(\Gamma\cdot S),\quad(x,g)\mapsto gf(g^{-1}x).
$$
Then $\phi$ is equivariant, and for $g,h\in\Gamma$ with $g^{-1}h\in E$ satisfies
\begin{align}\label{contract}
d(\phi(x,g),\phi(x,h)) &  =d(gf(g^{-1}x),hf(h^{-1}x)) \nonumber \\ & =d(f(g^{-1}x),g^{-1}hf(h^{-1}x)) \nonumber \\ & <\epsilon.
\end{align}
Now, let $\{V_0,...,V_d\}$ be the cover of $C$ given by Lemma \ref{nice cover}.  With notation from Lemma \ref{nice cover}, define
$$
\mathcal{U}:=\{\phi^{-1}(N_{\frac{1}{6}10^{-d}}(V_{i~\Delta}))~|~i\in \{0,...,d\} \text{ and } \Delta \text{ an $i$-simplex}\}.
$$
The open cover $\mathcal{U}$ of $X\times \Gamma$ is equivariant as it is the pullback of an equivariant cover by an equivariant map.  It satisfies condition (A) as $\Gamma$ permutes the disjoint sets $N_{\frac{1}{6}10^{-d}}(V_{i~\Delta})$ for each fixed $i$; it satisfies condition $(B)$ as $C$ is a $(\Gamma,\mathcal{F})$-complex; it satisfies condition (C) as each point in $C$ can intersect at most one $N_{\frac{1}{6}10^{-d}}(V_{i~\Delta})$ for each $i$; and it satisfies condition (D) as $S$ is finite.  Moreover, as $\{V_0,...,V_d\}$ covers $C$, the Lebesgue number of the cover 
$$
\mathcal{V}:=\{N_{\frac{1}{6}10^{-d}}(V_{i~\Delta})~|~i\in \{0,...,d\} \text{ and } \Delta \text{ an $i$-simplex}\}.
$$
of $C$ is at least $\epsilon=\frac{1}{6}10^{-d}$.  Hence the condition
$$
g^{-1}h\in E \ra d(\phi(x,g),\phi(x,h))<\epsilon
$$
from line \eqref{contract} above implies for any $x$ and $g$ there exists $V\in\mathcal{V}$ such that 
$$
\phi(\{x\}\times gE)\subseteq V;
$$
this implies condition (E) for $\mathcal{U}$.\\

Now assume condition (ii).   Let $\epsilon,E$ be given; replacing $E$ by $E\cup E^{-1}\cup\{e\}$, we may assume that $E$ is symmetric and contains the identity element.  Let $d$ be as in (ii), and take a cover $\mathcal{U}$ of $X\times\Gamma$ as in (ii) for the finite set $E^n=\{g\in \Gamma~|~g=g_1...g_n \text{ for some } g_1,...,g_n\in E\}$ where $n$ is such that
$$
\frac{(2d+2)(4d+6)}{n}<\epsilon.
$$ 
Let $C$ be the \emph{nerve} of this cover: precisely, $C$ is the subcomplex of $P(\mathcal{U})$ consisting of the union of the convex hulls of all subsets $\{U_0,...,U_m\}$ of the vertex set $\mathcal{U}$ such that $\bigcap_{i=0}^m U_i\neq \varnothing$.  Note that $C$ is a $(\Gamma,\mathcal{F})$-complex by conditions (A) and (B); it is moreover of dimension at most $d$ by (C).  

The next part of the argument is a `topological version' of \cite[Proposition 4.1]{Dadarlat:2007qy}.  For each $U\in \mathcal{U}$ and each $m\in \{0,...n\}$ define the \emph{$E^m$-interior of $U$} to be 
$$
U^{(m)}:=\{(x,g)\in U~|~\{x\}\times gE^m\subseteq U\}.
$$
This gives rise to a nested sequence of open sets 
$$
U^{(n)}\subseteq U^{(n-1)}\subseteq \cdots \subseteq U^{(1)}\subseteq U^{(0)}=U.
$$
Note that for each $m\in \{0,...,n\}$, $\{U^{(m)}\}_{U\in \mathcal{U}}$ is an open cover of $X\times\Gamma$ by condition (E).  Let $\{V_U\}_{U\in\mathcal{U}}$ be an open cover of $X\times \Gamma$ such that for each $U$, $\overline{V_U}\subseteq U^{(n)}$; as $X\times \Gamma$ is normal, standard arguments in general topology imply that such a `shrunken' cover exists. 

Fix $U\in\mathcal{U}$ for the moment.   Define $V^{(n)}:=V_U$.  Define 
$$
W^{(n)}:=\{(x,gh)\in X\times \Gamma~|~(x,g)\in V^{(n)},h\in E\}
$$
and note that if $(x,g)$ is in $\overline{W^{(n)}}$, then $(x,gh)$ is in $V^{(n)}$ for some $h\in E$, and thus 
$$
\{x\}\times gE^{n-1}\subseteq \{x\}\times ghE^n\subseteq  U.
$$ 
This says that $\overline{W^{(n)}}\subseteq U^{(n-1)}$, and thus by normality of $X\times \Gamma$ there is an open set $V^{(n-1)}$ with 
$$
\overline{W^{(n)}}\subseteq V^{(n-1)}\subseteq \overline{V^{(n-1)}}\subseteq U^{(n-1)}.
$$
Now set 
$$
W^{(n-1)}:=\{(x,gh)\in X\times \Gamma~|~(x,g)\in V^{(n-1)},h\in E\},
$$
and similarly we have $\overline{W^{(n-1)}}\subseteq U^{(n-2)}$ whence there is an open set $V^{(n-2)}$ with 
$$
\overline{W^{(n-1)}}\subseteq V^{(n-2)}\subseteq \overline{V^{(n-2)}}\subseteq U^{(n-2)}.
$$
Continuing this process gives a nested sequence 
$$
V^{(n)}\subseteq V^{(n-1)}\subseteq \cdots\subseteq V^{(1)}\subseteq V^{(0)}
$$
of open sets such that: for each $m\in \{1,...,n\}$, $\overline{V^{(m)}}\subseteq V^{(m-1)}$; for each $m\in \{0,...,n\}$, $\overline{V^{(m)}}\subseteq U^{(m)}$; and for each $m\in \{1,...,n\}$,
\begin{equation}\label{v alm inv}
\{(x,gh)\in X\times \Gamma~|~(x,g)\in V^{(m)},h\in E\}\subseteq V^{(m-1)}.
\end{equation}
Now, for each $m\in \{1,...,n\}$, Urysohn's lemma implies there exists a continuous function $\psi_U^{(m)}:X\times \Gamma\to [0,1]$ such that 
$$
\psi_U^{(m)}(V^{(m)})=\{1\}\quad \text{and}\quad \psi_U^{(m)}((X\times \Gamma)\setminus V^{(m-1)})=\{0\}.
$$
Define 
$$
\psi_U:=\sum_{m=1}^n \psi_U^{(m)}.
$$
Note that if $(x,g)\in V^{(m-1)}\setminus V^{(m)}$ for some $m\in \{1,...,n\}$ then $\psi_U(x,g)\in [m-1,m]$ and moreover $\psi_U(x,g)=0$ whenever $(x,g)\not\in V^{(0)}$.  Hence by line \eqref{v alm inv} we have that for any $(x,g)\in X\times \Gamma$
\begin{equation}\label{first bound}
h\in E\quad\Rightarrow \quad |\psi_U(x,g)-\psi_U(x,gh)|\leq 2.
\end{equation}

At this point, we unfix $U$, and define $\phi_U:X\times \Gamma\to [0,1]$ by the formula
$$
\phi_U=\frac{\psi_U}{\sum_{V\in \mathcal{U}}\psi_V}.
$$
As $\{V_U\}_{U\in \mathcal{U}}$ is an open cover of $X\times \Gamma$, $\{\phi_U\}_{U\in \mathcal{U}}$ is a well-defined partition of unity on $X\times \Gamma$.  Moreover, for any $x\in X$, any $g,h\in G$ such that $g^{-1}h\in E$, and any $U\in \mathcal{U}$ we have
\begin{align*}
| & \phi_U(x,g)-\phi_U(x,h)|\\ & \leq \Big|\frac{\psi_U(x,g)-\psi_U(x,h)}{\sum_{V\in \mathcal{U}}\psi_V(x,g)}\Big|+\frac{\psi_U(x,h)}{{\sum_{V\in \mathcal{U}}\psi_V(x,h)}}\Big|\frac{\sum_{V\in \mathcal{U}}\psi_V(x,g)-\sum_{V\in \mathcal{U}}\psi_V(x,h)}{\sum_{V\in \mathcal{U}}\psi_V(x,g)}\Big| \\
& \leq \Big|\frac{\psi_U(x,g)-\psi_U(x,h)}{\sum_{V\in \mathcal{U}}\psi_V(x,g)}\Big|+\Big|\frac{\sum_{V\in \mathcal{U}}\psi_V(x,g)-\psi_V(x,h)}{\sum_{V\in \mathcal{U}}\psi_V(x,g)}\Big|. 
\end{align*}
By construction, at least one of the terms in the sum $\sum_{V\in \mathcal{U}}\psi_V(x,g)$ equals $n$.  Combining this with line \eqref{first bound} and the fact that the cover has multiplicity at most $d+1$ implies that the above is bounded by
$$
\frac{2}{n}+\frac{4(d+1)}{n}=\frac{4d+6}{n}.
$$
Hence for $x\in X$ and $g,h\in\Gamma$ with $g^{-1}h\in E$
$$
\sum_{U\in\mathcal{U}}|\phi_U(x,g)-\phi_U(x,h)|<\frac{(2d+2)(4d+6)}{n}.
$$
as at most $2d+2$ of the terms can be non-zero by the multiplicity restriction.  

Finally, define
$$
f:X\to C,~~~x\mapsto \sum_{U\in\mathcal{U}}\phi_U(x,e)U,
$$
which is continuous as each $\phi_U$ is.  Note that for $g\in E$, $x\in X$,
$$
d(f(gx),gf(x))=\sum_{U\in\mathcal{U}}|\phi_U(gx,e)-\phi_U(gx,g)|<\frac{(2d+2)(4d+6)}{n}<\epsilon
$$
by choice of $n$; thus we may conclude (i).
\end{proof}

\begin{definition}\label{blrdef}
If an action $\Gamma\lefttorightarrow X$ and family $\mathcal{F}$ satisfy the conditions in Proposition \ref{blr}, we say that the action is \emph{$d$-BLR for $\mathcal{F}$}.
\end{definition}

The following theorem shows that the conditions on an action in Proposition \ref{blr} imply that the action has dynamic asymptotic dimension at most $d$.   We do not know whether the converse is true: some evidence for a converse is provided by Corollary \ref{fadcor} below that relates the $d$-BLR condition and dynamic asymptotic dimension to asymptotic dimension in the sense of coarse geometry.  

\begin{theorem}\label{blrlp}
Let $\Gamma\lefttorightarrow X$ be an action with $X$ compact, and let $\mathcal{F}$ be the family of finite subgroups of $\Gamma$.  If $\Gamma\lefttorightarrow X$ is a $d$-BLR action for $\mathcal{F}$, then it has dynamic asymptotic dimension at most $d$.
\end{theorem}

\begin{proof}
Assume condition (i) in Proposition \ref{blr} for the family $\mathcal{F}$ of finite subgroups of $\Gamma$.   Let a finite subset $E$ of $\Gamma$ be given; we may assume that $E$ contains the identity and that $E=E^{-1}$.  Let $\epsilon=\frac{1}{3}10^{-d}$.  Then there exists a $(\Gamma,\mathcal{F})$-complex of dimension at most $d$ and an $(E,\epsilon)$-equivariant continuous map
$$
f:X\to C
$$ 
as in Definition \ref{aemaps}.   Using Remark \ref{perturb}, we may assume that there exists a finite subset $S$ of the vertex set of $C$ such that $f(X)$ is  contained in $C\cap P(S)$.  Define also 
\begin{equation}\label{f def}
F:=\{g\in \Gamma~|~gS\cap S\neq \varnothing\}.
\end{equation}
Note that as $S$ is finite and as the stabilizer in $\Gamma$ of each vertex in $S$ is finite, $F$ is a finite subset of $\Gamma$.  

Now, let $\{V_0,...,V_d\}$ be the open cover of $C$ as in Lemma \ref{nice cover}; we will freely use the notation from Lemma \ref{nice cover} in the rest of the proof.  For each $i\in \{0,...d\}$ define $U_i:=f^{-1}(V_i)$ and for each $i$-simplex $\Delta$, define $U_{i~\Delta}:=f^{-1}(V_{i~\Delta})$.  Note that $\{U_0,...,U_d\}$ is an open cover of $X$, and that each $U_i$ is the disjoint union of the sets $U_{i~\Delta}$ as $\Delta$ ranges over all $i$-simplices.  We claim that the following holds for any $i$ and $i$-simplex $\Delta$:
\begin{equation}\label{key point}
x\in U_{i~\Delta},~gx\in U_i,~g\in E\quad \Rightarrow \quad gx\in U_{i~g\Delta}.
\end{equation}
Indeed, with notation as in the line above, we know that $gx$ is in $U_{i~\Delta'}$ for some $i$-simplex $\Delta'$.  As $f$ is $(E,\epsilon)$-equivariant we have
$$
d(f(gx),gf(x))<\epsilon;
$$
As $gf(x)$ is in $V_{i~g\Delta}$ and $f(gx)$ is in $V_{i~\Delta'}$, line \eqref{separation} from Lemma \ref{nice cover} now forces $g\Delta=\Delta'$ and thus $gx$ is in $U_{i~g\Delta}$ as claimed.

To complete the proof, say $x\in U_{i}$ for some $i$ and $g=g_n\cdots g_1$ with $g_k\in E$ and $g_k\cdots g_1x\in U_i$ for each $k\in \{1,...,n\}$.  If $x$ is in $U_{i~\Delta}$ for some $i$-simplex $\Delta$, then repeated applications of line \eqref{key point} force $gx$ to be in $U_{i~g\Delta}$, so in particular $f(gx)\in V_{i~g\Delta}\cap f(X)$.  As both $f(gx)$ and $f(x)$ are supported on vertices in $S$, this in turn forces $gS\cap S$ to be non-empty, and thus $g$ is in the finite set $F$ from line \eqref{f def} above.  This shows that $\Gamma\lefttorightarrow X$ satisfies the conditions in Definition \ref{daddef} so we are done.
\end{proof}

\begin{remark}\label{damen}
In earlier versions of this work, we considered the following definition.  An action $\Gamma\lefttorightarrow X$ with $X$ compact is \emph{$d$-amenable} if for all finite subsets $E$ of $\Gamma$ and all $\epsilon>0$ there exists an $(E,\epsilon)$-equivariant map 
$$
f:X\to P_d(\Gamma).
$$ 
This definition is motivated by the notion of an amenable action, which is the same definition without the extra restriction on $d$.  We also said the \emph{amenability dimension} of an action is the smallest $d$ such that the action is $d$-amenable.  This notion was used in recent work of Szab\'{o}, Wu, and Zacharias, where they gave very interesting relations to Rokhlin-type conditions for actions of residually finite groups: see \cite{Szabo:2014aa}, particularly Section 4.

Clearly a $d$-amenable action is $d$-BLR; it is not difficult to see that the converse is true for torsion free groups, but false (for example) for finite groups in general (take $X$ a point, and $\Gamma$ to a finite group of cardinality larger than $d$).  Thus in particular, Theorem \ref{blrlp} also implies that $d$-amenable actions have dynamic asymptotic dimension at most $d$ (and thus by our other results, $d$-amenability has consequences for nuclear dimension).  
\end{remark}

\section{Dynamic asymptotic dimension for groupoids}\label{gpdsec}

In this section we reformulate the definition of dynamic asymptotic dimension for groupoids.  Our main goal is to explore connections to coarse geometry in Section \ref{asdimsec}, but also the extra generality seemed interesting for the results in Section \ref{fndsec}.

For us, groupoids are always locally compact and Hausdorff; we generally leave these assumptions implicit from now on.  We will also assume that groupoids are \'{e}tale, but do not leave this (less standard) assumption implicit.   Accessible, as well as reasonably self-contained and concise, introductions to this class of groupoids and their associated $C^*$-algebras can be found in \cite[Section 2.3]{Renault:2009zr} and \cite[Section 5.6]{Brown:2008qy}.  

We use the following notation for groupoids.  A groupoid is denoted $G$, and its unit space by $G^{(0)}$.  We identify $G^{(0)}$ with a closed and open subspace of $G$ in the canonical way.  The range and source maps are denoted $r,s:G\to G^{(0)}$.  An ordered pair $(g,h)$ of elements from $G$ is \emph{composable} if $s(g)=r(h)$, and the composition is written $gh$ in this case.  The inverse of $g\in G$ is denoted $g^{-1}$.

\begin{definition}\label{gpddad}
Let $G$ be an \'{e}tale groupoid.  Then $G$ has \emph{dynamic asymptotic dimension} $d\in \N$ if $d$ is the smallest number with the following property: for every open relatively compact subset $K$ of $G$ there are open subsets $U_0,...,U_d$ of $G^{(0)}$ that cover $s(K)\cup r(K)$ such that for each $i$, the set $\{g\in K~|~s(g),r(g)\in U_i\}$ is contained in a relatively compact subgroupoid of $G$.
\end{definition}

Note that we may equivalently ask that for each $i$, the subgroupoid $G_i$ generated by the set $\{g\in K~|~s(g),r(g)\in U_i\}$ is relatively compact.  It is often convenient to note that this $G_i$ is automatically open.  This follows from the next quite general (and presumably well-known) lemma, which will be used several times in the remainder of the paper.

\begin{lemma}\label{open rel com}
Let $G$ be an \'{e}tale groupoid.  
\begin{enumerate}[(i)]
\item If $K$, $H$ are open \textup{(}respectively compact, respectively relatively compact\textup{)} subsets of $G$, then 
$$
K\cdot H:=\{kh\in G~|~k\in K,h\in H,~s(k)=r(h)\}
$$
is open \textup{(}respectively compact, respectively relatively compact\textup{)}.
\item If $K\subseteq G$ is open, then the subgroupoid of $G$ generated by $K$ is also open, and is itself an \'{e}tale groupoid.
\end{enumerate}
\end{lemma}

\begin{proof}
For part (i), assume first that $K$ and $H$ are open, and let $kh$ be a point of $KH$.  As $G$ is \'{e}tale there are neighborhoods $V\subseteq K$ of $k$ and $W\subseteq H$ of $h$ such that $r$ and $s$ both restrict to homeomorphisms on $V$ and $W$.  Define $V':=(s|_V)^{-1}(s(V)\cap r(W))$, $W':=(r|_W)^{-1}(s(V)\cap r(W))$, and  
$$
U:=\{k'h'~|~k\in V',h\in W', s(k')=r(h')\};
$$
note that $U$ is a subset of $KH$.  As $s$ and $r$ both restrict to homeomorphisms on $V'$ and $W'$, $r$ restricts to a homeomorphism from $U$ to $r(V')$.  In particular, $U$ contains an open neighborhood of $kh$, so we are done.  Assume next that $K$ and $H$ are compact.  It is clear from the corresponding property for $K$ and $H$ that any net in $K\cdot H$ has a convergent subnet, so $K\cdot H$ is compact.  Finally, assume that $K$ and $H$ are relatively compact.  Clearly $K\cdot H\subseteq \overline{K}\cdot \overline{H}$, and we have already shown that the latter set is compact.

For (ii), note that if $H$ is any open subset of $G$ then 
$$
H^n:=\{h_n\cdots h_1~|~h_k\in H \text{ and }s(h_k)=r(h_{k-1})\text { for all } k\in\{1,...,n\} \}
$$
is open by (i).  If $K$ is open, then so too is $K\cup K^{-1}\cup r(K)\cup s(K)$ and the subgroup generated by $K$ is equal as a set to
$$
\bigcup_{n=1}^\infty\big(K\cup K^{-1}\cup r(K)\cup s(K)\big)^n
$$
and thus open.  An open subgroupoid of an \'{e}tale groupoid is easily seen to be itself \'{e}tale, so this completes the proof.
\end{proof}

We give a first, easy, example: it is a groupoid analogue of Remark \ref{dadrem}, parts \eqref{dad0} and \eqref{dad02}.

\begin{example}\label{lfgpd}
A groupoid is \emph{locally finite} if it is the union of open, relatively compact subgroupoids.  It is not difficult to check directly from Definition \ref{gpddad} that a groupoid has dynamic asymptotic dimension zero if and only if it is locally finite in this sense. 

In particular, if $G=G^{(0)}$ is just a space, then $G$ has dynamic asymptotic dimension zero.
\end{example}

We now show that Definition \ref{gpddad} generalizes our earlier definition for group actions.  An action $\Gamma\lefttorightarrow X$ gives rise to an associated \emph{transformation groupoid} $\Gamma\ltimes X$ as follows.  To fix notation, recall that as a set $\Gamma\ltimes X$ is defined by
$$
\Gamma\ltimes X:=\{(gx,g,x)~|~x\in X,g\in\Gamma\}.
$$
Note that the projection onto the last two variables is a bijection from $\Gamma\ltimes X$ to $\Gamma\times X$; $\Gamma\ltimes X$ is given the topology that makes this bijection a homeomorphism.  The unit space of $G^{(0)}$ is $\{(x,e,x)~|~x\in X\}$ and the source and range maps are given by
$$
r(gx,g,x)=gx,~~~s(gx,g,x)=x.
$$
Composition and inverse are defined by
$$
(ghx,g,hx)(hx,h,x)=(ghx,gh,x),~~~(gx,g,x)^{-1}=(x,g^{-1},gx)
$$
respectively.

\begin{lemma}\label{gpgpd}
An action $\Gamma\lefttorightarrow X$ has dynamic asymptotic dimension $d$ in the sense of Definition \ref{daddef} if and only if the corresponding transformation groupoid $G=\Gamma\ltimes X$ has dynamic asymptotic dimension $d$ in the sense of Definition \ref{gpddad}.
\end{lemma}

\begin{proof}
Say an action $\Gamma\lefttorightarrow X$ has dynamic asymptotic dimension $d$, and let $K$ be an open relatively compact subset of the transformation groupoid $\Gamma\ltimes X$.  Any open relatively compact subset of $G$ is contained in a (relatively compact, open) set of the form $\{(gx,g,x)\in G~|~x\in K',g\in E\}$ for some finite subset $E$ of $\Gamma$ such that $e\in E$ and $E=E^{-1}$, and some open  relatively compact subset $K'$ of $X$.  Hence we may assume that $K$ has this form.  Let $U_0',...,U_d'$ be open subsets of $G^{(0)}$ that cover $\overline{K'}$ and are such that the set
$$
F_i:=\left\{\begin{array}{l|l}  & \text{there exist } x\in U_i'\text{ and }g_n,...,g_1\in E\text{ such that } \\ g\in \Gamma &g=g_n\cdots g_2g_1 \text{ and  for each } k\in \{1,...,m\}, \\ &  g_k\cdots g_1x\in U_i'\end{array}\right\}
$$
is finite for each $i$.  Define $U_i:=U_i'\cap K'$, so $U_1,...,U_d$ are open subsets of $G^{(0)}=X$ that cover $K'$.  Then the subgroupoid of $G$ generated by 
$$
\{(gx,g,x)\in K~|~x,gx,\in U_i\}
$$
is contained in $\{(gx,g,x)\in G~|~g\in F_i,x\in K'\}$, and is thus relatively compact.

The converse can be proved in a very similar way. 
\end{proof}

\section{Example - coarse spaces with finite asymptotic dimension}\label{asdimsec}

In this section we show that the coarse groupoid $G(X)$ associated to a bounded geometry coarse space $X$ has dynamic asymptotic dimension $d$ if and only if $X$ has asymptotic dimension $d$.  This example (in the special case $X$ is the coarse space underlying a group $\Gamma$) was our original motivation, and also motivates the terminology `dynamic asymptotic dimension'.

We recall the following definition: compare for example \cite[Chapter 2]{Roe:2003rw}.

\begin{definition}\label{cs}
Let $X$ be a set.  A \emph{coarse structure on $X$} is a collection $\mathcal{E}$ of subsets of $X\times X$ called \emph{controlled sets} such that:
\begin{enumerate}[(i)]
\item the diagonal $\{(x,x)~|~x\in X\}$ is contained in $\mathcal{E}$;
\item if $E$ is in $\mathcal{E}$ and $F$ is a subset of $E$, then $F$ is in $\mathcal{E}$;
\item if $E$ and $F$ are in $\mathcal{E}$, then their union $E\cup F$ is in $\mathcal{E}$;
\item if $E$ is in $\mathcal{E}$ then
$$
E^{-1}:=\{(x,y)\in X\times X~|~(y,x)\in E\}
$$
is in $\mathcal{E}$;
\item if $E$ and $F$ are in $\mathcal{E}$, then 
$$
E\circ F:=\{(x,y)\in X\times X~|~\text{there exists }z\in X \text{ such that } (x,z)\in E, (z,y)\in F\}
$$
is in $\mathcal{E}$.
\end{enumerate}
A set $X$ equipped with a coarse structure in called a \emph{coarse space}. 

A coarse space $X$ has \emph{bounded geometry} if for each controlled set $E$ the cardinals
$$
\sup_{x\in X}|\{y\in X~|~(x,y)\in E\}|,~~~\sup_{x\in X}|\{y\in X~|~(y,x)\in E\}|
$$
are finite.
\end{definition}

The reader should keep the following example of a coarse space in mind.  If $(X,d)$ is a discrete metric space, a \emph{controlled tube} is defined to be a set of the form
$$
E:=\{(x,y)\in X\times X~|~d(x,y)<r\}
$$ 
for some $r\geq 0$, and the controlled sets are defined to be those subsets $E$ of $X\times X$ that are contained in a controlled tube. The bounded geometry condition corresponds to the following assumption: for any $r>0$, there is a uniform bound on the cardinality of all balls of radius $r$ in $X$.

We recall the following definition, due originally to Gromov \cite[Section 1.E]{Gromov:1993tr} in the special case that the coarse structure comes from a metric as above.  See for example \cite[Chapter 9]{Roe:2003rw}, \cite[Chapter 2]{Nowak:2012aa}, or \cite[Section 12]{Bell:2008fk} for more information.

\begin{definition}\label{asdim0}
Let $X$ be a coarse space and $E$ a controlled set for $X$.  A cover $\mathcal{U}=\{U_i\}_{i\in I}$ of $X$ is:
\begin{enumerate}[(i)]
\item \emph{$E$-separated} if whenever $i\neq j$, $U_i\times U_j \cap E=\varnothing$;
\item \emph{$E$-bounded} if each set $U_i\times U_i$ is contained in $E$. 
\end{enumerate}
\end{definition}

\begin{definition}\label{asdim}
Let $X$ be a coarse space.  Then $X$ has \emph{asymptotic dimension} $d$ if $d$ is the smallest number with the following property: for any controlled set $E$ there exists a controlled set $F$ and a cover $\mathcal{U}=\{U_i\}_{i\in I}$ of $X$ such that $\mathcal{U}$ is $F$-bounded and such that $\mathcal{U}$ admits a decomposition 
$$\mathcal{U}=\mathcal{U}_0\sqcup \cdots \sqcup \mathcal{U}_d$$
 such that each $\mathcal{U}_i$ is $E$-separated.
\end{definition}

Now, let $X$ be a bounded geometry coarse space.  Skandalis, Tu, and Yu \cite{Skandalis:2002ng} (see also \cite[Chapter 10]{Roe:2003rw}) associate a groupoid to $X$ as follows.  For each controlled set $E$ for $X$, let $\overline{E}$ denote its closure in the product $\beta X\times \beta X$ of the Stone-\v{C}ech compactification of $X$ with itself.  Define
$$
G(X):=\bigcup_{E\in\mathcal{E}} \overline{E},
$$
and equip $G(X)$ with the weak topology coming from this union, i.e.\ a subset $U$ of $G(X)$ is open if and only if its intersection with each $\overline{E}$ is open in the natural compact topology on $\overline{E}$.  This is a locally compact Hausdorff topology for which each $\overline{E}$ is a compact open subset of $G(X)$.  It is not difficult to see that if $(\omega,\eta)$ are in $\overline{E}$, and $(\eta,\zeta)$ are in $\overline{F}$, then $(\omega,\zeta)$ is in $\overline{E\circ F}$, whence it follows that $G(X)$ inherits a groupoid structure from the pair groupoid structure on $\beta X\times \beta X$ (its topology, however, is \emph{not} the same as the one it inherits from $\beta X\times \beta X$).  The groupoid operations satisfy all the necessary continuity axioms to show that $G(X)$ is a (locally compact, Hausdorff)  \'{e}tale groupoid.

\begin{theorem}\label{fadthe}
A bounded geometry coarse space $X$ has asymptotic dimension $d$ if and only if the associated coarse groupoid $G(X)$ has dynamic asymptotic dimension $d$.

Moreover, if $X$ has asymptotic dimension at most $d$, then $G(X)$ has dynamic asymptotic dimension at most $d$ in the following slightly stronger form: for any open, relatively compact subset $K$ of $G(X)$, there exists an open cover $\{U_0,...,U_d\}$ of $G(X)^{(0)}$ such that the subset
$$
\{g\in K~|~s(g),r(g)\in U_i\}
$$
is contained in a compact open subgroupoid of $G(X)$.
\end{theorem}

\begin{proof}
Assume first that $G(X)$ has dynamic asymptotic dimension $d$.  Fix a controlled set $E$ as in the definition of asymptotic dimension, which we may assume contains the diagonal, and let $K=\overline{E}$, a compact open subset of $G(X)$ that contains the unit space $G(X)^{(0)}=\beta X$.  Take an open cover $U_0,...,U_d$ of $G(X)^{(0)}=\beta X$ with the properties in Definition \ref{gpddad} for this $K$, so for each $i\in \{0,...,d\}$ the subgroupoid $G_i$ generated by the set 
$$
\{g\in \overline{E}~|~s(g),r(g)\in U_i\}
$$
is relatively compact.  In particular, by definition of the topology on $G(X)$ there exists a controlled set $F\subseteq X\times X$ such that each $G_i$ is contained in $\overline{F}\subseteq G(X)$.  Fixing $i$ for the moment, let $\sim$ be the equivalence relation on $U_i$ induced by $G_i$, so $x\sim y$ if there exists $g\in G_i$ with $s(g)=x$ and $r(g)=y$.  Write the equivalence classes for this relation as $\{U_i^j~|~j\in J_i\}$ and define $\mathcal{U}_i:=\{U_i^j\cap X~|~j\in J_i\}$ to be the collection of equivalence classes for this relation.  Now let $\mathcal{U}=\mathcal{U}_0\sqcup \cdots \sqcup \mathcal{U}_d$ be the collection of all these subsets of $X$ (which covers $X$ as $\{U_0,...,U_d\}$ covers $\beta X$).

We claim that $\mathcal{U}$, decomposed as $\mathcal{U}=\mathcal{U}_0\sqcup \cdots \sqcup \mathcal{U}_d$, has the properties required by Definition \ref{asdim}.  Indeed, for $F$-boundedness, note that the equivalence relation induced on each $U_i^j\cap X$ by $G_i$ is entirely contained in $F$ (as $G_i$ is contained in $\overline{F}$).  For $E$-separatedness, note that if for $j\neq k$ there was some $(x,y)$ in $(U_i^j\times U_i^k)\cap E$, then we would have that $x\sim y$, contradicting the fact that $U_i^j$ and $U_i^k$ are distinct equivalence classes.

Conversely, say $G$ has asymptotic dimension at most $d$.  Let $K$ be a compact subset of $G(X)$.  The definition of the topology on $G(X)$ implies that $K\subseteq \overline{E}$ for some controlled set $E$.   Let $\mathcal{U}=\mathcal{U}_0\sqcup \cdots \sqcup \mathcal{U}_d$ and $F$ be as in the definition of finite asymptotic dimension for this $E$.  For each $i$, set
$$
U_i:=\overline{\bigsqcup_{U\in \mathcal{U}_i}U},
$$
which is a compact open subset of $\beta X$; as $\mathcal{U}$ is a cover of $X$, $\{U_0,...,U_d\}$ is a cover of $\beta X$.  Note that for each $i$, 
$$
G_{i}:=\bigsqcup_{U\in \mathcal{U}_i}U\times U
$$
is a subgroupoid of the pair groupoid that is contained (as a set) in $F$, and by continuity of the groupoid operations, the subgroupoid of $G(X)$ generated by 
$$
\{g\in E~|~r(g),s(g)\in U_i\}
$$ 
is contained in $\overline{G_i}$, a subgroupoid of $G(X)$ contained (as a set) in the compact subset $\overline{F}$.  Hence the dynamic asymptotic dimension of $G(X)$ is at most $d$.

The final statement is clear from the construction in the proof.
\end{proof}

Finally in this section, we give two consequences.  The first discusses the relationship of the above to the Bartels-L\"{u}ck-Reich conditions of Section \ref{blrsec}; the second shows that groups of finite asymptotic dimension all admit free, amenable, minimal actions on the Cantor set which have finite dynamic asymptotic dimension.

\begin{theorem}\label{fadcor}
Let $\Gamma$ be a discrete group.  The following are equivalent:
\begin{enumerate}[(i)]
\item $\Gamma$ admits an action $\Gamma\lefttorightarrow X$ on a compact space that is $d$-BLR with respect to the family of finite subgroups.
\item $\Gamma$ admits an action $\Gamma\lefttorightarrow X$ on a compact space with dynamic asymptotic dimension at most $d$.
\item The canonical action of $\Gamma$ on $\beta\Gamma$ has dynamic asymptotic dimension at most $d$.
\item $\Gamma$ equipped with the canonical left-invariant coarse structure for which the controlled sets are
$$
\mathcal{E}=\{E\subseteq \Gamma\times \Gamma~|~\{s^{-1}t\in \Gamma~|~(s,t)\in E\} \text{ is finite}\},
$$
has asymptotic dimension at most $d$.
\item The canonical action of $\Gamma$ on $\beta\Gamma$ is $d$-BLR with respect to the family of finite subgroups.
\end{enumerate}
\end{theorem}

\begin{proof}
The fact that (i) implies (ii) (for the same $X$) is Theorem \ref{blrlp}.  

For (ii) implies (iii), assume that $\Gamma\lefttorightarrow X$ is an action on a compact space with dynamic asymptotic dimension at most $d$.  Let $\Gamma\to X$ be any orbit map, and let $\phi:\beta\Gamma\to X$ be the canonical extension to the Stone-\v{C}ech compactification of $\Gamma$.  Pulling back covers along $\phi$ shows that $\Gamma\lefttorightarrow\beta \Gamma$ also has dynamic asymptotic dimension at most $d$.

For (iii) implies (iv), note that the transformation groupoid $\Gamma\ltimes \beta\Gamma$ has dynamic asymptotic dimension at most $d$ by Lemma \ref{gpgpd}.  On the other hand, this transformation groupoid canonically identifies with the coarse groupoid of $\Gamma$ by \cite[Proposition 3.4]{Skandalis:2002ng}, so the result follows from Theorem \ref{fadthe}.

The fact that (iv) implies (v) is essentially due to Higson and Roe, and very similar to the proof of Proposition \ref{blr}.  The result follows from combining the arguments of \cite[Lemma 4.3]{Higson:2000dp} and  \cite[Section 3]{Higson:2000dp}: indeed, these give an $(E,\epsilon)$-equivariant map $\beta\Gamma\to P_d(\Gamma)$ for any finite $E\subseteq \Gamma$ and any $\epsilon>0$.  

Finally, (v) implies (i) is trivial.
\end{proof}

The following theorem is a result of combining Theorem \ref{fadthe} with ideas of R\o{}rdam and Sierakowski from \cite[Section 6]{Rordam:2010kx}.

\begin{theorem}\label{mincor}
Let $\Gamma$ be a countably infinite discrete group with asymptotic dimension $d$.  Then $\Gamma$ admits a free, minimal action on the Cantor set which has dynamic asymptotic dimension at most $d$.
\end{theorem}

\begin{proof}
Using Theorem \ref{fadcor}, the action of $\Gamma$ on $\beta\Gamma$ is $d$-BLR for the family $\mathcal{F}$ of finite subgroups of $\Gamma$.  Let 
$$
E_1\subseteq E_2\subseteq E_3\subseteq \cdots
$$
be a nested sequence of finite subsets of $\Gamma$ such that $\Gamma=\bigcup E_n$.  For each $n$, let $f_n:\beta\Gamma\to C_n$ be a $(E_n,1/n)$-equivariant map to a $(\Gamma,\mathcal{F})$ simplicial complex of dimension at most $d$.  As $\beta\Gamma$ is compact, we may assume by Remark \ref{perturb} that the image of each $f_n$ is contained in a set of the form $P(S_n)\cap C_n$, where $S_n$ is a finite subset of the vertex set of $C_n$.  Hence, replacing $C_n$ by $C_n\cap P(\Gamma\cdot S_n)$, we may assume that there are only finitely many $\Gamma$-orbits of vertices of $C_n$, and thus in particular that the vertex set $V_n$ of $C_n$ is countable.

Now, let $A$ denote a unital $\Gamma$-invariant $C^*$-subalgebra of $l^\infty(\Gamma)$.  Let $\widehat{A}$ denote the spectrum of $A$, which by Gelfand duality identifies with a quotient space of $\beta\Gamma$ such that the quotient map $\beta\Gamma\to \widehat{A}$ is equivariant. For each $n$, write 
$$
f_n(x)=\sum_{v\in V_n}t_v(x)v,
$$
where each $t_v:\beta \Gamma\to [0,1]$ is continuous.  Define $\Omega_1:=\{t_v~|~v\in V_n,n\in \N\}$, which is a countable subset of $l^\infty(\Gamma)$.  If $A$ contains $\Omega_1$, then the action $\Gamma\lefttorightarrow \widehat{A}$ is $d$-BLR: indeed, the functions $t_v$ then descend to the quotient space $\widehat{A}$ of $\beta\Gamma$, whence the maps $f_n:\beta\Gamma\to C_n$ descend to $f_n:\widehat{A}\to C_n$  (this argument is inspired by \cite[Lemma 3.5]{Higson:2000bl}).   On the other hand, \cite[Lemma 6.4]{Rordam:2010kx} shows that there is a countable subset $\Omega_2$ of $l^\infty(\Gamma)$ such that if $A$ contains $\Omega_2$, then the action of $\Gamma$ on $\widehat{A}$ is free.

Define $\Omega:=\Omega_1\cup \Omega_2$, a countable subset of $l^\infty(\Gamma)$.  R\o{}rdam and Sierakowski \cite[Lemma 6.7]{Rordam:2010kx} show there is a countable $\Gamma$-invariant collection of projections $\mathcal{P}$ such that the $C^*$-algebra generated by $\mathcal{P}$ contains $\Omega$.  Let $A$ be the $C^*$-algebra generated by $\mathcal{P}$ and the unit of $l^\infty(G)$.  It follows from the above discussion that $\widehat{A}$ is a $\Gamma$-space, and that the action of $\Gamma$ is free and $d$-BLR.  Moreover, as $A$ is unital and generated by countably many projections, the space $\widehat{A}$ is compact, metrizable, and totally disconnected.

We now have that $\Gamma$ admits a free, $d$-BLR action on a totally disconnected, metrizable compact space $Y:=\widehat{A}$.  Note that if $Z$ is any non-empty closed $\Gamma$-invariant subset of $Y$, then $Z$ and the induced $\Gamma$ action on it will still have all these properties.  Inclusion defines a partial order on the closed non-empty $\Gamma$-invariant subsets of $Y$; compactness implies that any descending chain for this order has non-empty intersection, and thus Zorn's lemma implies there exists a minimal element $X$.  It follows easily that the induced $\Gamma$ action on $X$ is minimal.  We now have that $\Gamma$ admits a free, minimal, $d$-BLR action on a totally disconnected, metrizable compact space $X$.  

To complete the proof, it suffices to show that this $X$ is a copy of the Cantor set, and for this it suffices to show that it has no isolated points.  Note then that if $x\in X$ is isolated then the orbit of $\Gamma\cdot x$ is open, and thus $X\setminus \Gamma\cdot x$ is empty by minimality.  As the action is free and $\Gamma$ is infinite, the open cover $\{\{gx\}~|~g\in \Gamma\}$ of $X=\Gamma\cdot x$ has no finite subcover which contradicts compactness.  Hence no isolated point can exist, and we are done.
\end{proof}

\section{Partitions of unity}\label{pousec}

In this section we prove a technical result showing that dynamic asymptotic dimension gives rise to partitions of unity which are `almost' invariant in an appropriate sense.  This will be important in Section \ref{fndsec}.

\begin{proposition}\label{dadpou}
Let $G$ be an \'{e}tale groupoid with compact unit space, and with dynamic asymptotic dimension $d$.  Then for any open relatively compact subset $K$ of $G$ and any $\epsilon>0$ there exists an open cover $\{U_0,...,U_d\}$ of $r(K)\cup s(K)$ with the following properties.
\begin{enumerate}[(i)]
\item \label{generate0} For each $i$, the set
$$
\{g\in K~|~s(g),r(g)\in U_i\}
$$
is contained in an \textup{(}open and\textup{)} relatively compact subgroupoid of $G$.
\item \label{fat} For each $x\in r(K)\cup s(K)$, the `partial orbit' $s(r^{-1}(x)\cap K)$ is completely contained in some $U_i$.
\item \label{pou} There exists a collection of continuous functions $\{\phi_i:G^{(0)}\to [0,1]\}_{i=0}^d$ such that the support of each $\phi_i$ is contained in $U_i$, such that for all $x\in r(K)\cup s(K)$ we have $\sum_{i=0}^d (\phi_i(x))^2=1$, and such that for any $g\in K$ and each $i$
$$
|\phi_i(s(g))-\phi_i(r(g))|<\epsilon.
$$
\end{enumerate}
\end{proposition}

For the convenience of readers who are mainly interested in the case of group actions on compact spaces, we spell out here what this says in the special case where $G$ is a transformation groupoid in the language of actions.  The proposition is not really any easier to prove in this special case, however.

\begin{corollary}\label{dadgppou}
Let $\Gamma\lefttorightarrow X$ be an action with $X$ compact, and with dynamic asymptotic dimension $d$.  Then for any finite subset $E$ of $\Gamma$ and any $\epsilon>0$ there exists an open cover $\{U_0,...,U_d\}$ of $X$ with the following properties.
\begin{enumerate}[(i)]
\item For each $i$, the set
$$
\left\{\begin{array}{l|l}  g\in \Gamma & \exists x\in U_i \text{ and } g_n,...,g_1\in E\text{ such that } g=g_n\cdots g_1 \\ & \text{and } \forall k\in \{1,...,n\},~g_k\cdots g_1x\in U_i\end{array}\right\}
$$
is finite.
\item For each $x\in X$, the collection $E\cdot x:=\{g^{-1}x~|~g\in E\}$ is completely contained in some $U_i$.
\item There exists a collection of continuous functions $\{\phi_i:X\to [0,1]\}_{i=0}^d$ on $X$ such that the support of each $\phi_i$ is contained in $U_i$, such that $\sum_{i=0}^d (\phi_i(x))^2=1$ for all $x\in X$, and such that for any $g\in E$ and each $i$
$$
\sup_{x\in X}|\phi_i(gx)-\phi_i(x)|<\epsilon. \eqno\qed
$$
\end{enumerate}
\end{corollary}

The proof of Proposition \ref{dadpou} is related to that of Proposition \ref{blr}.  We wanted to keep the two proofs separate for the convenience of readers who are only interested in one or the other case, and as there are enough significant differences that it did not seem possible to make a `combined' proof that was much shorter.

We start with a lemma.

\begin{lemma}\label{gpddad2}
Say $G$ is an \'{e}tale groupoid with dynamic asymptotic dimension $d$.  Then for any open relatively compact subset $K$ of $G$ there is a cover $\{U_0,...,U_d\}$ of $r(K)\cup s(K)$ by relatively compact open subsets of $G^{(0)}$ such that for each $i$, the set $\{g\in K~|~s(g),r(g)\in U_i\}$ generates an \textup{(}open and\textup{)} relatively compact subgroupoid of $G$, and moreover so that for each $x\in G^{(0)}$ the set $s(r^{-1}(x)\cap K)$
is completely contained in some $U_i$.  
\end{lemma}

\begin{proof}
Replacing $K$ with $K\cup K^{-1}\cup r(K)\cup s(K)$ (which is open and relatively compact as $r$ and $s$ are open, continuous maps), we may assume that $K=K\cup K^{-1}\cup r(K)\cup s(K)$.  Define 
$$
K^3:=\{g\in G~|~\text{ there are } g_1,g_2,g_3\in K \text{ such that } g=g_1g_2g_3\},
$$
which is again an open relatively compact subset of $G$ by Lemma \ref{open rel com}.  Let $\{V_0,...,V_d\}$ be an open cover of $r(K^3)\cup s(K^3)$ such that the subgroupoid $G_i$ of $G$ generated by 
\begin{equation}\label{subgpd}
\{g\in K^3~|~s(g),r(g)\in V_i\}
\end{equation}
is relatively compact.   For $i\in \{0,...,d\}$, define
$$
U_i:=s(K\cap r^{-1}(V_i))\cap (r(K)\cup s(K)),
$$
which is open and relatively compact.  Note that as $r(K)\cup s(K)\subseteq K$, each $U_i$ contains $V_i\cap ( r(K)\cup s(K) )$; thus we have 
\begin{align*}
\bigcup_{i=0}^d U_i &  \supseteq \bigcup_{i=0}^d V_i\cap ( r(K)\cup s(K) )\supseteq (r(K^3)\cup s(K^3))\cap (r(K)\cup s(K))  \\ & =r(K)\cup s(K).
\end{align*}

It remains to show that for each fixed $i$, the set $\{g\in K~|~s(g),r(g)\in U_i\}$ generates a relatively compact subgroupoid of $G$.  Indeed, say $g$ is in the subgroupoid generated by this set.  Then (recalling that $K=K\cup K^{-1}\cup r(K)\cup s(K)$) there exists a finite sequence $g_n,...,g_1$ of composable elements from $\{g\in K~|~s(g),r(g)\in U_i\}$ such that $g=g_n\cdots g_1$.  As for each $k\in \{1,...,n\}$ we have that $s(g_k)\in U_i$, there exists $h_k\in K$ such that $r(h_k)\in V_i$, and $s(h_k)=s(g_k)$.  Moreover, as $s(g_n)\in U_i$, there exists $h_{n+1}\in K$ such that $s(h_{n+1})=s(g_n)$ and such that $r(h_{n+1})\in V_i$.  For each $k\in \{1,...,n\}$, define
$$
g_k':=h_{k+1}g_kh_k^{-1},
$$
and note that $g_k'\in \{g\in K^3~|~s(g),r(g)\in V_i\}$.  On the other hand, $g=h_{n+1}^{-1}g_n'\cdots g_1'h_1$.  If $G_i$ is the subgroupoid of $G$ generated by the set in line \eqref{subgpd}, it follows that $g$ is in the subset 
$$
K\cdot G_i \cdot K :=\{h\in G~|~\text{there are } h_0,h_2\in K \text{ and }h_1\in G_i \text { such that } h=h_0h_1h_2\};
$$
by relative-compactness of $G_i$ and $K$, and Lemma \ref{open rel com}, this set is relatively compact, so we are done.
\end{proof}

We need one more technical lemma before giving the proof of Proposition \ref{dadpou}; in some sense it is an iterated version of Lemma \ref{gpddad2} above.

\begin{lemma}\label{gen top}
Say $G$ is an \'{e}tale groupoid with dynamic asymptotic dimension $d$.  Let $K$ be an open relatively compact subset of $G$, and let $N$ be a fixed natural number.   Then for each $i\in \{0,...d\}$ there is a nested collection 
$$
U_i^{(0)}\subseteq U_i^{(1)}\subseteq \cdots \subseteq U_i^{(N+1)}
$$ 
of open, relatively compact subsets of $G^{(0)}$ with the following properties.
\begin{enumerate}[(i)]
\item\label{cover2} The collection $\{U_0^{(0)},...,U_d^{(0)}\}$ covers $\overline{r(K)\cup s(K)}$.
\item\label{closure} For all $i,n$, $\overline{U_i^{(n)}}\subseteq U_i^{(n+1)}$.
\item\label{bigger} For all $i,n$, 
$$
s(K\cap r^{-1}(U_i^{(n)}))\subseteq U_i^{(n+1)}.
$$
\item\label{generate2} For all $i$, the set
$$
\{g\in K~|~s(g),r(g)\in U_i^{(N+1)}\}
$$
generates a relatively compact subgroupoid of $G$.
\end{enumerate}
\end{lemma}

\begin{proof}
Replacing $K$ with $K\cup K^{-1}\cup r(K) \cup s(K)$, we may assume that $K=K\cup K^{-1}\cup r(K) \cup s(K)$.  For each $n$, define
$$
\overline{K}^n:=\{g_n\cdots g_1~|~g_k\in \overline{K} \text{ and } s(g_{k+1})=r(g_k) \text{ for all k}\}.
$$
Note that $\overline{K}^n$ is compact for all $n$ by Lemma \ref{open rel com}, and that $\overline{K}^n\subseteq \overline{K}^{n+1}$ for all $n$ by our assumption that $K=K\cup K^{-1}\cup r(K) \cup s(K)$.  As $G$ is locally compact, there exists an open, relatively compact subset $K'\subseteq G$ that contains $\overline{K}^{N+1}$.  Let $\{V_0,...,V_d\}$ be an open cover of $r(K')\cup s(K')$ with the properties in Lemma 7.3 for the open relatively compact set $K'$.  Replacing each $V_i$ with its intersection with $r(K')\cup s(K')$, we may assume that each $V_i$ is contained in $r(K')\cup s(K')$.

Now, for each $i\in \{0,...d\}$ and each $n\in \{0,...,N+1\}$, define 
$$
V_i^{(n)}:=\{x\in V_i~|~s(\overline{K}^{N+1-n}\cap r^{-1}(x))\subseteq V_i\}.
$$
As the sets $\overline{K}^n$ are nested, this gives a nested sequence of subsets 
$$
V_i^{(0)}\subseteq V_i^{(1)}\subseteq \cdots \subseteq V_i^{(N+1)}
$$
of $G^{(0)}$.  As for each $x\in G^{(0)}$, $s(r^{-1}(x)\cap K')\subseteq V_i$ for some $i$, the collection $\{V_0^{(0)},...V_d^{(0)}\}$ covers $r(K')\cup s(K')$.  We claim moreover that each $V_i^{(n)}$ is open.  Indeed, if this does not happen, then we may find a net $(x_\lambda)$ in $G^{(0)}$ converging to some $x\in V_i^{(n)}$ such that for each $\lambda$ there exists $g_\lambda\in  \overline{K}^{N+1-n}\cap r^{-1}(x_\lambda)$ with $s(g_\lambda)\not\in V_i$.  As $\overline{K}^{N+1-n}$ is compact we may assume by passing to a subnet that $(g_\lambda)$ converges to some $g\in \overline{K}^{N+1-n}$.  As $r$ is continuous, $g$ is in $ \overline{K}^{N+1-n}\cap r^{-1}(x)$, and as $V_i$ is open, $s(g)$ is not in $V_i$.  This contradicts that $x$ is in $V_i^{(n)}$, so $V_i^{(n)}$ is open as claimed.

We may now define the sets $U_i^{(n)}$.  As the set $\overline{r(K')\cup s(K')}$ is compact, it is normal.  It follows that there are open subsets $U_0^{(0)},...,U_d^{(0)}$ of $r(K')\cup s(K')$ such that $\overline{U_i^{(0)}}\subseteq V_i^{(0)}$ for each $i$ and such that the collection $\{U_0^{(0)},...,U_d^{(0)}\}$ covers $\overline{r(K)\cup s(K)}$.  Note that for each $i$,
\begin{align*}
\overline{s(K\cap r^{-1}(U_i^{(0)}))} & \subseteq s(\overline{K\cap r^{-1}(U_i^{(0)})})\subseteq s(\overline{K}\cap \overline{r^{-1}(U_i^{(0)})})\subseteq s(\overline{K}\cap r^{-1}(\overline{U_i^{(0)}})) \\ & \subseteq s(\overline{K}\cap r^{-1}(V_i^{(0)})).
\end{align*}
It is not difficult to see that this last set is contained in $V_i^{(1)}$, however.  Moreover, $\overline{U_i^{(0)}}$ is contained in $V_i^{(0)}\subseteq V_i^{(1)}$ by assumption.  Hence by normality of the compact set $\overline{V_i}$, there exists an open set $U_i^{(1)}$ containing $\overline{s(K\cap r^{-1}(U_i^{(0)}))}\cup \overline{U_i^{(0)}}$ such that $\overline {U_i^{(1)}}\subseteq V_i^{(1)}$.  Continuing in this way, we see that 
$$
\overline{s(K\cap r^{-1}(U_i^{(1)}))} \subseteq V_i^{(2)},
$$
and thus there is an open set $U_i^{(2)}$ containing $\overline{s(K\cap r^{-1}(U_i^{(1)}))}\cup \overline{U_i^{(1)}}$, with closure contained in $V_i^{(2)}$, and so on.  

The process above gives a nested sequence of open subsets
$$
U_i^{(0)}\subseteq U_i^{(1)}\subseteq \cdots \subseteq U_i^{(N+1)}
$$
for each $i$ such that 
$$
\overline{s(K\cap r^{-1}(U_i^{(n)}))} \subseteq U_i^{(n+1)} \quad \text{and} \quad \overline{U_i^{(n)}}\subseteq V_i^{(n)}
$$
for all $n$.  Note that each $U_i^{(n)}$ is relatively compact as $V_i^{(n)}\subseteq V_i\subseteq r(K')\cup s(K')$ for all $n$.  It remains to check properties (i) through (iv) from the statement. Properties \eqref{cover2}, \eqref{closure}, and \eqref{bigger} are obvious from the way we chose our sets.  Property \eqref{generate2} follows as
$$
\{g\in K~|~s(g),r(g)\in U_i^{(N+1)}\}\subseteq \{g\in K'~|~s(g),r(g)\in V_i\},
$$
and the set on the right generates a relatively compact subgroupoid by choice of the sets $V_0,...,V_d$. 
\end{proof}

\begin{proof}[Proof of Proposition \ref{dadpou}]
Let $K$ be an open, relatively compact subset of $G$ and $\epsilon>0$ be given.  Replacing $K$ with $K\cup K^{-1}\cup r(K) \cup s(K)$, we may assume that $K=K^{-1}$, and that $K$ contains its image under $r$ and $s$.  Let $N$ be any natural number larger than $2$ such that $(\sqrt{2}(1+\sqrt{d+1}))/\sqrt{N}<\epsilon$, and let $\{U_i^{(n)}~|~i\in \{0,...,d\},~n\in \{0,...,N+1\}\}$ have the properties in Lemma \ref{gen top}.  For each $i=0,...,d$, set $U_i=U_i^{(N+1)}$.  Condition \eqref{cover2} from Lemma \ref{gen top} combined with the assumptions from the first paragraph of that lemma imply that $\{U_0,...U_d\}$ is an open cover of $r(K)\cup s(K)$.  Condition \eqref{fat} from the statement of Proposition \ref{dadpou} follows condition \eqref{bigger} from Lemma \ref{gen top} and the fact that $\{U_0^{(0)},...,U^{(0)}_d\}$ covers $r(K)\cup s(K)$.  Condition \eqref{generate0} from the statement of Proposition \ref{dadpou} follows directly from condition \eqref{generate2} from Lemma \ref{gen top}.  It remains to construct a partition of unity with the properties in condition \eqref{pou} from the statement of Proposition \ref{dadpou}.

For each $i\in\{0,...,d\}$ and $n\in \{1,...,N\}$, let
$$
\psi_i^{(n)}:G^{(0)}\to [0,1]
$$ 
be any continuous function which is constantly equal to one on $U_i^{(n-1)}$ and constantly equal to zero on $G^{(0)}\setminus U_i^{(n)}$; such functions exist by condition \eqref{closure} from Lemma \ref{gen top} and Urysohn's lemma.  For $i\in \{0,...,d\}$ set 
$$
\psi_i=\frac{1}{N}\sum_{n=1}^N \psi_i^{(n)}
$$
and define
$$
\phi_i=\frac{\psi_i}{\max\Big\{\sqrt{\sum_{j=0}^d \psi_j^2},1\Big\}}.
$$
Note that condition \eqref{cover2} from Lemma \ref{gen top} implies that for all $x\in r(K)\cup s(K)$, there is at least one $i\in \{0,...,d\}$ such that $\psi_i(x)=1$, whence 
$$
\sum_{i=0}^d (\phi_i(x))^2=1
$$
for all $x\in r(K)\cup s(K)$.

Fix now $g\in K$ and $i\in \{0,...,d\}$; to complete the proof we must show that 
\begin{equation}\label{lip}
|\phi_i(s(g))-\phi_i(r(g))|<\epsilon.
\end{equation}
For notational convenience, for each $j\in \{0,...,d\}$, set $U_j^{(n)}=G^{(0)}$ for $n\geq N+2$.  Define
$$
M=M_j:=\min\{n~|~r(g)\in U_j^{(n)}\}.
$$
Note that 
$$
\psi_j^{(n)}(r(g))=\left\{\begin{array}{ll} 1 & n\geq M+1 \\ 0 & n\leq M-1\end{array}\right.
$$
and $0\leq\psi(r(g))\leq 1$ whence
\begin{equation}\label{range}
\frac{N-(M+1)}{N}\leq \psi_j(r(g))\leq\frac{N-M}{N}.
\end{equation}
Note also that condition \eqref{bigger} from Lemma \ref{gen top} combined with the fact that $K=K^{-1}$ implies that $s(g)$ is in $U_j^{(M+1)}\setminus U_j^{(M-2)}$ whence 
$$
\psi_j^{(n)}(s(g))=\left\{\begin{array}{ll} 1 & n\geq M+2 \\ 0 & n\leq M-2\end{array}\right.;
$$
as $0\leq \psi_j^{(n)}(s(g))\leq 1$ for all values of $n$, this forces
\begin{equation}\label{source}
\frac{N-(M+2)}{N}\leq \psi_j(s(g))\leq \frac{N-(M-1)}{N}.
\end{equation}
Combining lines \eqref{range} and \eqref{source}, we may conclude that 
\begin{equation}\label{psibound}
|\psi_j(r(g))-\psi_j(s(g))|\leq \frac{2}{N}
\end{equation}

One the other hand, at least one of the $\psi_j$ is equal to one on each of $s(g)$ and $r(g)$, and therefore
\begin{equation}\label{below}
1\leq \sum_{j=0}^d\psi_j(r(g)) ~~\text{ and }~~1\leq \sum_{j=0}^d\psi_j(s(g)).
\end{equation}
Hence for our fixed choice of $i\in \{0,...,d\}$ and $g\in K$,
\begin{align*}
&|\phi_i(r(g))-\phi_i(s(g))| =\Bigg|\frac{\psi_i(r(g))}{\sqrt{\sum_{j=0}^d \psi_j^2(r(g))}}-\frac{\psi_i(s(g))}{\sqrt{\sum_{j=0}^d \psi_j^2(s(g))}}\Bigg| 
\end{align*}
and so
\begin{align}\label{main bound}
|\phi_i(r(g))-\phi_i(s(g))| &  \leq \frac{1}{\sqrt{\sum_{j=0}^d \psi_j^2(r(g))}} |\psi_i(r(g))-\psi_i(s(g))|\nonumber \\ & \quad+|\psi_i(s(g))|\Bigg|\frac{1}{\sqrt{\sum_{j=0}^d \psi_j^2(r(g))}}-\frac{1}{\sqrt{\sum_{j=0}^d \psi_j^2(s(g))}}\Bigg|.
\end{align}
Using lines \eqref{psibound} and \eqref{below} we have 
\begin{equation}\label{first term}
\frac{1}{\sqrt{\sum_{j=0}^d \psi_j^2(r(g))}} |\psi_i(r(g))-\psi_i(s(g))| \leq \frac{2}{N}.
\end{equation}
On the other hand,
\begin{align*}
& |\psi_i(s(g))|\Bigg|\frac{1}{\sqrt{\sum_{j=0}^d \psi_j^2(r(g))}}-\frac{1}{\sqrt{\sum_{j=0}^d \psi_j^2(s(g))}}\Bigg| \\
&=  \frac{1}{\sqrt{\sum_{j=0}^d \psi_j^2(r(g))}}\frac{1}{\sqrt{\sum_{j=0}^d \psi_j^2(s(g))}}\Bigg|\sqrt{\sum_{j=0}^d \psi_j^2(s(g))}-\sqrt{\sum_{j=0}^d \psi_j^2(r(g))}\Bigg| \\
& \leq \Bigg|\sqrt{\sum_{j=0}^d \psi_j^2(s(g))}-\sqrt{\sum_{j=0}^d \psi_j^2(r(g))}\Bigg| \\
& \leq \sqrt{\Big|\sum_{j=0}^d \psi_j^2(s(g))-\sum_{j=0}^d \psi_j^2(r(g))\Big|} \\
& \leq \sqrt{\sum_{j=0}^d |\psi_j(s(g))-\psi_j(r(g))||\psi_j(s(g))+\psi_j(r(g))|}.
\end{align*}
Using line \eqref{psibound} again thus gives
\begin{equation}\label{second term}
|\psi_i(s(g))|\Bigg|\frac{1}{\sqrt{\sum_{j=0}^d \psi_j^2(r(g))}}-\frac{1}{\sqrt{\sum_{j=0}^d \psi_j^2(s(g))}}\Bigg| \leq \sqrt{\frac{4(d+1)}{N}}.
\end{equation}
Combining lines \eqref{main bound}, \eqref{first term} and \eqref{second term} therefore gives
\begin{align*}
|\phi_i(r(g))-\phi_i(s(g))| & \leq \frac{2}{N}+\sqrt{\frac{4(d+1)}{N}}\leq \sqrt{\frac{2}{N}}+\sqrt{\frac{4(d+1)}{N}} \\ & =\frac{\sqrt{2}(1+\sqrt{d+1})}{\sqrt{N}}.
\end{align*}
which is smaller than $\epsilon$ by choice of $N$.
\end{proof}

\section{Nuclear dimension}\label{fndsec}

As usual, we adopt the conventions that all groupoids are locally compact and Hausdorff.  Recall that a groupoid $G$ is \emph{free} (such groupoids are also called \emph{principal}) if for each $x\in G^{(0)}$, the \emph{isotropy group} defined by
$$
G_x^x:=\{g\in G~|~s(g)=r(g)=x\}
$$
is trivial.  

In this section, we study nuclear dimension for the reduced groupoid $C^*$-algebra $C^*_r(G)$ of a free groupoid $G$ of finite dynamic asymptotic dimension.  Modulo a minor technical assumption that is satisfied for example when $G$ is second countable, our main result---Theorem \ref{fnd} below---says that the nuclear dimension of $C^*_r(G)$ is bounded above by $(d+1)(N+1)-1$, where $d$ is the dynamic asymptotic dimension of $G$, and $N$ is the covering dimension of $G^{(0)}$.   The result is inspired by (and implies) the Winter-Zacharias result that the nuclear dimension of the uniform Roe algebra associated to a bounded geometry metric space is at most the asymptotic dimension of the space \cite[Theorem 8.5]{Winter:2010eb}.  See the end of the section for some other corollaries.

We first recall the definitions of nuclear dimension, covering dimension, and of the reduced $C^*$-algebra of a groupoid.  More details on nuclear dimension can be found in the paper of Winter and Zacharias \cite{Winter:2010eb} that introduces the notion.  We do not really need a definition of covering dimension, but record it as there is some ambiguity about which definition one should use for non-metrizable spaces.  Accessible introductions to groupoid $C^*$-algebras, that largely focus on the case of interest here, can be found in \cite[Section 2.3]{Renault:2009zr} and  \cite[Section 5.6]{Brown:2008qy}.

\begin{definition}\label{nddef}
A completely positive map $\phi:A\to B$ is  \emph{order zero} if it preserves orthogonality: in other words, if $a_1,a_2$ are positive elements of $A$ such that $a_1a_2=0$, then $\phi(a_1)\phi(a_2)=0$.

Let $A$ be a $C^*$-algebra.  The \emph{nuclear dimension} of $A$ is the smallest integer $d\in \N$ with the following property.  For any finite subset $\mathcal{F}$ of $A$ and any $\epsilon>0$ there exist finite dimensional $C^*$-algebras $F_0,...,F_d$ and contractive completely positive maps
$$
\xymatrix{ A \ar[dr]^{\Phi_i} & & A \\ & F_i \ar[ur]^{\Psi_i} & } 
$$
such that each $\Psi_i$ is order zero, and such that for all $a\in \mathcal{F}$,
$$
\Big\|\sum_{i=0}^d(\Psi_i\circ \Phi_i)(a)-a\Big\|<\epsilon\|a\|.
$$
\end{definition}

\begin{definition}\label{covdimdef}
Let $X$ be a paracompact (Hausdorff) topological space.  The \emph{covering dimension} of $X$ is the smallest integer $d\in \N$ such that every open cover of $X$ has a refinement $\mathcal{V}$ such that $\mathcal{V}$ splits into a disjoint union
$$
\mathcal{V}=\mathcal{V}_0\sqcup \cdots \sqcup \mathcal{V}_d
$$
such that whenever $V\neq W$ are distinct elements of some $\mathcal{V}_i$, then $V\cap W=\varnothing$.
\end{definition}

The definition above is equivalent to several other commonly used definitions of covering dimension on the class of paracompact Hausdorff spaces: compare for example \cite[Proposition 1.6]{Kirchberg:2004uq} or \cite[Remark 4.5]{Pears:1975ty}.  

\begin{definition}\label{gpdalg}
Let $G$ be an \'{e}tale groupoid.  Let $C_c(G)$ denote the vector space of compactly supported complex-valued functions on $G$, made into a $*$-algebra via the convolution product and adjoint defined by
$$
(f_1f_2)(g)=\sum_{g_1g_2=g}f_1(g_1)f_2(g_2),~~~f^*(g)=\overline{f(g^{-1})}.
$$
For each $x\in G^{(0)}$ let $l^2(s^{-1}(x))$ denote the Hilbert space of square-summable functions on the source fibre of $x$, and define a $*$-representation $\pi_x$ of $C_c(G)$ on $l^2(s^{-1}(x))$ by
$$
(\pi_x(f)\xi)(g)=\sum_{g_1g_2=g}f(g_1)\xi(g_2).
$$
The \emph{reduced groupoid $C^*$-algebra of $G$} is the completion of $C_c(G)$ for the norm 
$$
\|f\|:=\sup_{x\in G^{(0)}}\|\pi_x(f)\|.
$$
\end{definition}

Before we state the main theorem, we need one (ad-hoc) definition.

\begin{definition}\label{small}
Let $G$ be an \'{e}tale groupoid, and $H$ be an open subgroupoid of $G$.  The subgroupoid $H$ is \emph{small} if it is either compact, or second countable and relatively compact.

A groupoid $G$ has \emph{strong dynamic asymptotic dimension at most $d$} if for any open relatively compact subset $K$ of $G$ there exists an open cover $\{U_0,...,U_d\}$ of $s(K)\cup r(K)$ such that for each $i$, the set
$$
\{g\in K~|~s(g),r(g)\in U_i\}
$$
is contained in a small subgroupoid of $G$.
\end{definition}

\begin{remark}\label{str rem}
\begin{enumerate}[(i)]
\item A second countable groupoid with dynamic asymptotic dimension at most $d$ automatically has strong dynamic asymptotic dimension at most $d$.  Thus readers who are only interested in the second countable case can just ignore the word `strong' in all the statements below, and replace `small' with `open and relatively compact'.
\item \label{fad str} Theorem \ref{fadthe} implies that if $X$ is a bounded geometry metric space of asymptotic dimension $d$, then the coarse groupoid $G(X)$ has strong dynamic asymptotic dimension $d$.
\item\label{dadpou still good} Say $G$ has strong dynamic asymptotic dimension $d$.  Then Proposition \ref{dadpou}, point \eqref{generate0} can be strengthened to say that $\{g\in K~|~s(g),r(g)\in U_i\}$ is contained in a small subgroupoid of $G$ for each $i$: indeed, exactly the same proof gives the stronger statement.
\end{enumerate}
\end{remark}

\begin{theorem}\label{fnd}
Let $G$ be a free, \'{e}tale, groupoid.  Assume that $G$ has strong dynamic asymptotic dimension at most $d$ \textup{(}Definition \ref{small}\textup{)}, and moreover that the unit space $G^{(0)}$ has topological covering dimension at most $N$.  

Then the nuclear dimension of the reduced groupoid $C^*$-algebra $C^*_r(G)$ is at most $(N+1)(d+1)-1$.
\end{theorem}

As well as the results subsumed by this theorem discussed in the corollaries at the end of this section, the reader might compare it to \cite[Theorem 4.6]{Szabo:2014aa}, which deduces analogous estimates on nuclear dimension from Rokhlin type conditions, and a condition related to asymptotic dimension.  We do not know the extent of the overlap between Theorem \ref{fnd} above and \cite[Theorem 4.6]{Szabo:2014aa}; the relationship between the hypotheses of the two results seems worth investigating more carefully.

Most of the rest of this section is devoted to the proof of Theorem \ref{fnd}; we give some corollaries at the end.  There are two main steps to the proof: an analysis of small subgroupoids, and a reduction to these subgroupoids.

\subsection*{Small subgroupoids}

The goal of this subsection is to prove the following fact about small subgroupoids.

\begin{proposition}\label{smallprop}
Let $G$ be a free, \'{e}tale groupoid with unit space.  Let $H$ be a small \textup{(}open\textup{)} subgroupoid of $G$ in the sense of Definition \ref{small}, and assume that its unit space $H^{(0)}$ has covering dimension $N$.  Then $C^*_r(H)$ has nuclear dimension\footnote{As will be clear from the proof, one can replace nuclear dimension with decomposition rank here, but we will not need this distinction.} $N$.
\end{proposition}

We will prove this separately in the case where $H$ is compact, and in the case where $H$ is relatively compact and second countable.  Before this point though, we need some basic definitions and lemmas that will be used in both cases.

\begin{definition}\label{equiv}
Let $H$ be a \textup{(}locally compact, Hausdorff\textup{)} groupoid.  Let $\sim$ be the equivalence relation on $H^{(0)}$ induced by $H$, i.e.\ $x\sim y$ if there is $h\in H$ with $s(h)=x$ and $r(h)=y$.  For $x\in H^{(0)}$, let $[x]$ denote its equivalence class.  For each $m\in \N$, define
$$
H^{(0)}_m:=\{x\in H^{(0)}~|~|[x]|=m\}.
$$
Finally, let $H^{(0)}/H$ and $H^{(0)}_m/H$ denote the spaces of equivalence classes equipped with the quotient topology in inherited from $H^{(0)}$.
\end{definition}

\begin{lemma}\label{quot open}
Let $H$ be a topological groupoid.  Then the quotient map $\pi:H^{(0)}\to H^{(0)}/H$ is open.
\end{lemma}

\begin{proof}
Say $U$ is an open subset of $H^{(0)}$.  Then as $s$ is continuous and $r$ is open, the set $r(s^{-1}(U))$ is open in $H^{(0)}$.  However, it is clear that $\pi^{-1}(\pi(U))=r(s^{-1}(U))$, and thus by definition of the quotient topology, $\pi(U)$ is open.
\end{proof}

\begin{lemma}\label{finite}
Let $G$ be an \'{e}tale groupoid and $H$ be a relatively compact subgroupoid of $G$.  Then there exists $M\in \N$ such that $H^{(0)}_m$ is empty for all $m> M$. 
\end{lemma}

\begin{proof}
Note that as $G$ is \'{e}tale, each $g\in G$ is contained in an open neighbourhood $U$ on which both $r$ and $s$ are injective.  As the closure $\overline{H}$ is compact, there are finitely many of these neighbourhoods covering $\overline{H}$, say $U_1,...,U_M$.  For each $x\in H^{(0)}$ we have $|[x]|\leq |(s|_H)^{-1}(x)|$; as $(s|_H)^{-1}(x)$ can intersect each $U_i$ at most once, however, $|(s|_H)^{-1}(x)|$ is bounded above by $M$.
\end{proof}

We will need the following classical theorem from dimension theory: see \cite[Proposition 2.16]{Pears:1975ty}.

\begin{proposition}\label{quotdim}
Let $X$ and $Y$ be paracompact Hausdorff topological spaces.  Let $\pi:X\to Y$ be a continuous open surjection such that for all $y\in Y$, $\pi^{-1}(y)$ is finite.  Then the covering dimensions of $X$ and $Y$ are equal. \qed
\end{proposition}

At this point, we specialize to the case of Proposition \ref{smallprop} where $H$ is compact; we will come back to the second countable and relatively compact case later.

\begin{lemma}\label{comdimcor}
Let $H$ be a compact \'{e}tale groupoid, and assume that the covering dimension of $H^{(0}$ is $N$.  Then $H^{(0)}/H$ is compact and Hausdorff, and the covering dimension of $H^{(0)}/H$ is exactly $N$.
\end{lemma}

\begin{proof}
As $H$ is compact and $H^{(0)}$ is a closed subspace of $H$, $H^{(0)}$ is compact.  Hence $H^{(0)}/H$ is a quotient space of a compact space, so compact.  To see that $H^{(0)}/H$ is Hausdorff, note that the equivalence relation $\sim$ on $H^{(0)}$ induced by $H$ is equal to 
$$
(r\times s)(H)\subseteq H^{(0)}\times H^{(0)},
$$
and is thus compact and in particular closed as $H^{(0)}\times H^{(0)}$ is Hausdorff.  It is a standard fact that the quotient of a compact space by an equivalence relation that is closed in this sense is Hausdorff: see for example \cite[Proposition 2.1]{Vick:1973pd}. The claim on covering dimension now follows immediately from Proposition \ref{quotdim} applied to the quotient map $\pi:H^{(0)}\to H^{(0)}/H$.
\end{proof}

\begin{proof}[Proof of Proposition \ref{smallprop} when $H$ is compact]
As $H$ is compact and $G$ is free, $H$ is a free and proper groupoid.  Hence $C^*_r(H)$ is Morita equivalent to $C(H^{(0)}/H)$ by \cite[Example 2.5 and Theorem 2.8]{Muhly:1987fk}\footnote{The cited paper only covers the second countable case, but the second countability assumption is unnecessary when the groupoid is \'{e}tale: see \cite{Felix:2014aa}.}.  Hence the nuclear dimension of $C^*_r(H)$ is equal to the covering dimension of $H^{(0)}/H$ by \cite[Proposition 2.4 and Corollary 2.8]{Winter:2010eb}, and this is $N$ by Lemma \ref{comdimcor}.
\end{proof}

We now turn to the case of Proposition \ref{smallprop} when $H$ is second countable with compact closure.  We start by recalling a result of Winter, which needs the following notation.  If $A$ is a $C^*$-algebra then $\text{Prim}(A)$ denotes the collection of all kernels of irreducible representations of $A$, equipped with the hull-kernel topology: this is defined by saying that a subset $S$ of $\text{Prim(A)}$ is closed if there exists an ideal $I$ of $A$ such that $S=\{J\in \text{Prim}(A)~|~J\supseteq I\}$.  Recall moreover that for $m\in \N$, $\text{Prim}_m(A)$ denotes the subspace of $\text{Prim}(A)$ consisting of kernels of $m$-dimensional representations; by \cite[Proposition 3.6.4]{Dixmier:1977vl}, $\text{Prim}_m(A)$ is locally compact and Hausdorff when equipped with the subspace topology.

Here then is a (weak version of) Winter's result: see \cite[Theorem 1.6]{Winter:2004rb}.

\begin{theorem}\label{shdr}
Let $A$ be a separable $C^*$-algebra such that all irreducible representations of $A$ have dimension at most $M$ for some $M\in \N$, and assume that $N$ is the maximal covering dimension of the spaces $\text{Prim}_m(A)$ for $m\in \{1,...,M\}$.  Then the nuclear dimension of $A$ is $N$. 
\end{theorem}

\begin{proof}
Winter shows in \cite[Theorem 1.6]{Winter:2004rb} that the decomposition rank of $C^*_r(H)$ is exactly $N$ under the stated assumptions.  Decomposition rank is (trivially - see \cite[Remarks 2.2 (ii)]{Winter:2010eb}) an upper bound for the nuclear dimension, and it follows from \cite[Proposition 2.9 and Corollary 2.10]{Winter:2010eb} that the nuclear dimension of $C^*_r(H)$ is bounded below by the maximum of the covering dimensions of the spaces $\text{Prim}_m(C^*_r(H))$.
\end{proof}

Clearly the remaining case of Proposition \ref{smallprop} follows from this theorem of Winter and the following result.

\begin{lemma}\label{scprops}
Let $H$ be a relatively compact, open subgroupoid of an \'{e}tale groupoid, and assume that the covering dimension of $H^{(0)}$ is $N$.  Then $\text{Prim}(C^*_r(H))$ is naturally homeomorphic to $H^{(0)}/H$, via a homeomorphism that takes $\text{Prim}_m(A)$ to $H^{(0)}_m/H$.  Moreover, the spaces $H^{(0)}_m/H$ are locally compact and Hausdorff, and the maximum of their covering dimensions is $N$.
\end{lemma}

For the proof of the lemma, we need a little more notation, and some preliminary lemmas about ideals in groupoid $C^*$-algebras.  Much of the following material seems likely to be standard for experts in groupoid $C^*$-algebras: for example, the results of \cite[Pages 101-103]{Renault:1980fk} and \cite[Corollary 4.9]{Renault:1991df} are closely connected to what follows.  However, we could not find complete proofs of exactly what we need in the literature, so give direct arguments below.

\begin{definition}\label{invop}
A subset $U$ of $H^{(0)}$ is \emph{invariant} if whenever $x\in U$ and $x\sim y$, we have $y\in U$.

We write $\mathcal{O}(H^{(0)})$ for the collection of all invariant open sets in $H^{(0)}$, and $\mathcal{I}(C^*_r(H))$ for the collection of all ideals in $C^*_r(H)$.  Both of these sets are equipped with the partial orders defined by inclusion.
\end{definition}

\begin{lemma}\label{lattice}
Let $H$ be an \'{e}tale groupoid.  For a subset $S$ of $C^*_r(H)$, write $\langle S\rangle$ for the ideal in $C^*_r(H)$ generated by $S$.  Provisionally define maps by
$$
\Phi:\mathcal{I}(C^*_r(H))\to \mathcal{O}(H^{(0)}),~~~I \mapsto \text{Prim}(I\cap C_0(H^{(0)}))
$$
and 
$$
\Psi:\mathcal{O}(H^{(0)})\to \mathcal{I}(C^*_r(H)),~~~U\mapsto \langle C_0(U)\rangle.
$$
Then $\Phi$ and $\Psi$ are well-defined and $\Phi(\Psi(U))=U$ for all $U\in \mathcal{O}(H^{(0)})$.
\end{lemma}

\begin{proof}
It is clear that $\Psi$ is well-defined.  To see that $\Phi$ is well-defined, we must show that the primitive ideal space $I\cap C_0(H^{(0)})$, which canonically identifies with the open subset 
$$
\{x\in H^{(0)}~|~\text{there exists } f\in I\cap C_0(H^{(0)})\text{ with } f(x)\neq 0\}
$$
of $H^{(0)}$, is invariant.  Say then $x$ is in this subset, and $x\sim y$, so there is $h\in H$ with $r(h)=x$ and $s(h)=y$.  Let $f\in I\cap C_0(H^{(0)})$ be such that $f(x)\neq 0$.  Let $U\owns h$ be an open subset of $H$ such that the restrictions of $r$ and $s$ to $U$ are injective.  Let $\phi:U\to [0,1]$ be any compactly supported function such that $\phi(h)=1$ (so in particular, $\phi$ is an element of $C^*_r(H)$).  Then $\phi f\phi^*$ is in $C_0(H^{(0)})\cap I$, and is non-zero on $y$; thus $y$ is in the primitive ideal space of $C_0(H^{(0)})\cap I$ as required.

We now check that for $U\in \mathcal{O}(H^{(0)})$ that $\Phi(\Psi(U))=U$.  This is equivalent to showing that 
$$
\langle C_0(U)\rangle \cap C_0(H)=C_0(U);
$$
as the inclusion $\langle C_0(U)\rangle \cap C_0(H)\supseteq C_0(U)$ is obvious, it remains to show the converse inclusion.  Say then $f$ is an element of $\langle C_0(U)\rangle \cap C_0(H^{(0)})$ and $x\in H^{(0)}$ is such that $f(x)\neq 0$; we will show that $x$ is in $U$.  As $f$ is in $\langle C_0(U)\rangle$ we may approximate it by elements of the algebraic ideal generated by $C_0(U)$, and thus in particular there must exist $f_i\in C_c(U)$ and $\phi_i,\psi_i\in C_c(H)$ for $i\in \{1,...,n\}$ such that 
$$
\sum_{i=1}^n (\phi_if_i\psi_i)(x)\neq 0.
$$
Hence for some fixed $i$, we have $(\phi_if_i\psi_i)(x)\neq 0$.  This says that 
$$
\sum_{g,h,k\in H,~x=ghk}\phi_i(g)f_i(h)\psi_i(k)\neq 0.
$$
As $f_i$ is supported in $U$, this implies that there must be $h\in U$ and $k\in H$ with $h=r(k)$ and $x=s(k)$.  Hence $x\sim h$, and thus $x$ is in $U$ by invariance.
\end{proof}

\begin{lemma}\label{lattice2}
Let $H$ be an open relatively compact subgroupoid of a free \'{e}tale groupoid $G$.  With notation as in Lemma \ref{lattice}, we have $\Psi(\Phi(I))=I$ for all $I\in \mathcal{I}(C^*_r(H))$.
\end{lemma}

\begin{proof} 
Given an ideal $I$ in $C^*_r(H)$, we need to show that $\langle I\cap C_0(H^{(0)})\rangle=I$.  The inclusion $I\supseteq \langle I\cap C_0(H^{(0)})\rangle$ is obvious, so it remains to show that $\langle I\cap C_0(H^{(0)})\rangle\supseteq I$.  Say then $f$ is an element of $I$.  Let $\epsilon>0$ and $f_0$ be an element of $C_c(H)$ such that $\|f-f_0\|_{C^*_r(H)}<\epsilon$ (note that $f_0$ need not be an element of $I$); write $K\subseteq H$ for the support of $f_0$.  

Fix for the moment $x\in H^{(0)}$ and write $s^{-1}(x)=\{h_1,...,h_m\}$ (this set is finite by Lemma \ref{finite}).  For each $h_i$ choose an open neighborhood $V_i$ of $h_i$ on which $r$ and $s$ restrict to homeomorphisms and which is such that $V_i\cap V_j=\varnothing$ for all $i\neq j$.  As $r(h_i)\neq r(h_j)$ for $i\neq j$, we may further assume that the sets $r(V_1),...,r(V_m)$ are mutually disjoint.  Set $V_x=\bigcap_{i=1}^m s(V_i)$, an open neighborhood of $x$.  

We now have an open cover $\{V_x\}_{x\in H^{(0)}}$ of $H^{(0)}$, and in particular of the compact subset $s(K)$.  Hence by standard results about existence of partitions of unity there is a finite collection $\{\phi_i:H^{(0)}\to [0,1]~|~i\in \{1,...,n\}\}$ of continuous functions such that each $\phi_i$ is supported in some compact subset of some $V_x$, and such that $\sum_{i=1}^n \phi_i(x)=1$ for all $x\in s(K)$.  It follows that $f_0(\sum_{i=1}^n \phi_i)=f_0$ (where the product is the convolution on $C_c(H)$).  As the norm of $C^*_r(H)$ restricts to the supremum norm on on the $C^*$-subalgebra $C_0(H^{(0)})$ this implies that
$$
\Big\|f\Big(\sum_{i=1}^n \phi_i\Big)-f\Big\|_{C^*_r(H)}\leq \|f-f_0\|_{C^*_r(H)}+\Big\|\sum_{i=1}^n \phi_i\Big\|_{C^*_r(H)}\|f-f_0\|_{C^*_r(H)}<2\epsilon.
$$
As $\epsilon$ was arbitrary, to complete the proof it suffices to show that each element $f\phi_i$ of $I$ is actually in $\langle I\cap C_0(H^{(0)})\rangle$.  

Fix then $\phi=\phi_i$, which is supported in a compact subset $K_x$ of some $V_x$.  Say $s^{-1}(x)=\{h_1,...,h_m\}$, and write $s^{-1}(V_x)=\bigsqcup_{j=1}^m V_j$, where $V_j$ is an open neighborhood of $h_j$ such that $s$ and $r$ restrict to homeomorphisms on each $V_j$, and so that the sets $r(V_1),...,r(V_m)$ are mutually disjoint (the existence of such $V_j$ follows from the construction of $V_x$).  For each $j$, let $\psi_j:H\to [0,1]$ be a continuous function supported in $V_j$, and such that $\psi_j(s^{-1}(K_x)\cap V_j)=\{1\}$ (such a $\psi_j$ exists by Urysohn's lemma).  Then $\sum_{j=1}^m\psi_j^*f\phi$ is an element of $C_0(H^{(0)})\cap I$.  Moreover, $f\phi=\sum_{j=1}^m\psi_j\psi_j^*f\phi$, so we are done.
\end{proof}

\begin{lemma}\label{prim}
Say $H$ is an open relatively compact subgroup of an \'{e}tale groupoid $G$.  Then the correspondence from Lemmas \ref{lattice} and \ref{lattice2} restricts to a bijection between the collection of primitive ideals in $C^*_r(H)$ and the collection of subsets of $H^{(0)}$ of the form $H^{(0)}\setminus [x]$ for some $x\in H^{(0)}$.
\end{lemma}

\begin{proof}
Note that by Lemma \ref{finite}, any equivalence class $[x]$ of some $x\in H^{(0)}$ is finite, so closed.  Hence all maximal open invariant sets are of the form $H^{(0)}\setminus [x]$ for some $x\in H^{(0)}$.  Hence by Lemmas \ref{lattice} and \ref{lattice2}, the maximal ideals in $\mathcal{I}(C^*_r(H))$ are exactly those of the form $\Psi(H^{(0)}\setminus [x])$ for some $x\in H^{(0)}$.  Maximal ideals are primitive by \cite[Theorem 2.9.7 (ii)]{Dixmier:1977vl}, so to complete the proof, it suffices to prove that any non-maximal ideal in $C^*_r(H)$ is not primitive.

We first claim that if $[x]$, $[y]$ are distinct equivalence classes in $H^{(0)}$, then there exist $U,V\in \mathcal{O}(H^{(0)})$ such that $[x]\subseteq U$, $[y]\subseteq V$, and such that for all $z\in U\cap V$, $|[z]|\geq |[x]|+|[y]|$.  Indeed, let $M$ be as in Lemma \ref{finite}.  Let $s:H\to H^{(0)}$ be the source map for $H$ (not for the ambient groupoid $G$), and write $s^{-1}(x)=\{h_1,...,h_m\}$ and $s^{-1}(y)=\{g_1,...,g_n\}$, where $|[x]|=m$, $|[y]|=n$ (this is possible as the groupoid $H$ is free).  As $[x]\cap[y]=\varnothing$ and $H$ is Hausdorff there are open sets $U_i\owns g_i$ and $V_j\owns h_j$ such that $s$ restricts to a homeomorphism on each $U_i$ and $V_j$, such that $U_i\cap U_j=\varnothing=V_i\cap V_j$ for $i\neq j$, and such that $U_i\cap V_j=\varnothing$ for all $i,j$.  Define
$$
U_0:=\bigcap_{i=1}^m s(U_i)~~\text{ and }~~V_0:=\bigcap_{j=1}^n s(V_j)
$$
and set
$$
U:=r(s^{-1}(U_0))~~\text{ and } ~~V:=r(s^{-1}(V_0)).
$$
Clearly $U$ and $V$ are open and invariant, and $[x]\subseteq U$, $[y]\subseteq V$.  Consider now $z\in U\cap V$.  Then $s^{-1}(z)$ must intersect all the sets $U_i$ and $V_j$; as these sets are disjoint, this forces $m+n\leq |s^{-1}(z)|$, and as $H$ is a free groupoid, this forces 
$$
|[z]|=|s^{-1}(z)|\geq  m+n=|[x]|+|[y]|,
$$ 
completing the proof of the claim.  

Now, say $V\in \mathcal{O}(H^{(0)})$ is not maximal, so there are $[x]\neq [y]$ with $[x],[y]\subseteq H^{(0)}\setminus V$.  Say without loss of generality $|[x]|\geq |[y]|\geq 1$.  The claim above implies there exist $U_x, U_y\in \mathcal{O}(H^{(0)})$ such that $[x]\subseteq U_x$ and $[y]\subseteq U_y$, and so that for any $z\in U_x\cap U_y$, $|[z]|>|[x]|$.  If $U_x\cap U_y=\varnothing$, it is not difficult to see that the images of the ideals $\Psi(U_x)$ and $\Psi(U_y)$ are non-zero and orthogonal in $C^*_r(H)/\Psi(V)$ whence $\Psi(V)$ is not primitive\footnote{If $A$ is a $C^*$-algebra faithfully represented on a Hilbert space $H$, and $I,J$ are non-zero orthogonal ideals in $A$, then $I\cdot H$ and $J\cdot H$ are $A$-invariant non-zero subspaces of $H$; in particular, the representation is reducible.} as required.  On the other hand, if there is some $z\in U_x\cap U_y$, then we may repeat the process with $[z]$ replacing $[x]$ and $[x]$ replacing $[y]$.  Using Lemma \ref{finite}, this process must finish eventually to give orthogonal non-zero ideals in $C^*_r(H)/\Psi(V)$, which is thus not primitive.  
\end{proof}

\begin{proof}[Proof of Lemma \ref{scprops}]
It follows from Lemma \ref{prim} that the map 
$$
H^{(0)}/H\to \text{Prim(A)},~~~[x]\mapsto \Psi(H^{(0)}\setminus [x])
$$
is a bijection.  Lemmas \ref{lattice} and \ref{lattice2}, the definition of the quotient topology on $H^{(0)}/H$, and the definition of the hull-kernel topology on $\text{Prim}(A)$ imply that this map is a homeomorphism.  Moreover, it is not difficult to see that $C^*_r(H)/\Psi(H^{(0)}\setminus [x])$ is isomorphic to the $C^*$-algebra $M_{|[x]|}(\C)$ of $|[x]|\times |[x]|$ matrices over $\C$, whence it follows that this homeomorphism takes $H_m^{(0)}/H$ onto $\text{Prim}_m(A)$.

Finally, note that the spaces $\text{Prim}_m(C^*_r(H))$ are locally compact and Hausdorff (\cite[Proposition 3.6.4]{Dixmier:1977vl}), and that for each $m$, $\text{Prim}_m(C^*_r(H))$ is open in $\sqcup_{n\leq m}\text{Prim}_n(C^*_r(H))$ (\cite[Proposition 3.6.3]{Dixmier:1977vl}).  Hence in particular $H^{(0)}_m$ is open in $\sqcup_{n\leq m} H^{(0)}_n$ (one could also prove this directly, of course).   As covering dimension does not increase on taking open subsets in our context (this follows for example from \cite[Proposition 2.5]{Winter:2010eb}), it follows inductively that the covering dimension of $H^{(0)}_m$ is at most $N$ for each $m$; on the other hand, it follows from \cite[Proposition 5.2]{Pears:1975ty} (plus second countability) that the covering dimension of $H^{(0)}$ is at most the maximum of the covering dimensions of $H^{(0)}_m$ for $m\in \{1,...,M\}$.  Hence the maximum of the covering dimensions of $H^{(0)}_m$, $m\in \{1,...,M\}$ is exactly $N$.  The claim on the covering dimension of $H^{(0)}_m/H$ now follows from Proposition \ref{quotdim} applied to the quotient map $\pi:H^{(0)}\to H^{(0)}/H$, which is open by Lemma \ref{quot open}.
\end{proof}

\subsection*{Completion of the proof of Theorem \ref{fnd}}

For this subsection, $G$ is as in the assumptions of Theorem \ref{fnd}.

For the next two lemmas, if $\phi$ and $f$ are in $C_c(G)$, let $f\cdot \phi$ denote their pointwise product in $C_c(G)$, i.e.\ $(\phi\cdot f)(g)=\phi(g)f(g)$.  Also, for a subset $K$ of $G$, we will write $C^*_K(G)$ for the subspace of $C^*_r(G)$ consisting of all elements supported in $K$.  

\begin{lemma}\label{multlem}
Let $K$ be a compact subset of $G$.  Then there exists a constant $M>0$ such that for all $\phi\in C^*_K(G)$ we have 
$$
\|\phi\cdot f\|_{C^*_r(G)}\leq M\sup_{g\in G}|\phi(g)|~\|f\|_{C^*_r(G)}.
$$
\end{lemma}

\begin{proof}
With notation as in Definition \ref{gpdalg}, it suffices to prove that there exists $M>0$ such that for all $f\in C_c(G)$ and all $x\in G^{(0)}$, 
$$
\|\pi_x(\phi\cdot f)\|\leq M\sup_{g\in G}|\phi(g)|~\|\pi_x(f)\|.
$$
As $G$ is \'{e}tale, for each $g\in G$ there is an open neighbourhood $U$ of $G$ such that both $r$ and $s$ are injective when restricted to $U$; as $K$ is compact, there is a collection $U_1,...,U_M$ of open subsets of $G$ with this property such that $K\subseteq \bigcup_{i=1}^M U_i$.  We may write $\phi$ as a sum $\phi=\phi_1+\cdots +\phi_M$, where each $\phi_i$ is supported in some $U_i$ as above, and satisfies $\sup_{g\in G}|\phi_i(g)|\leq \sup_{g\in G}|\phi(g)|$.  It thus suffices to prove that if $\psi\in C_c(G)$ is supported in an open set $U$ such that $r,s$ are injective when restricted to $U$, then
$$
\|\pi_x(\psi\cdot f)\|\leq \sup_{g\in G}|\psi(g)|\|\pi_x(f)\|.
$$

We now prove this.  Indeed, computing for $\xi\in l^2(s^{-1}(x))$ 
$$
\|\pi_x(\psi\cdot f)\xi\|^2=\sum_{g\in s^{-1}(x)}\Big|\sum_{g_1g_2=g}\psi(g_1)f(g_1)\xi(g_2)\Big|^2.
$$
For each $g\in G$, write $g_U$ for the unique element in $r^{-1}(r(g))\cap U$ (if this exists).  Then the sum above becomes
\begin{align*}
& \sum_{\{g\in s^{-1}(x)~|~r^{-1}(r(g))\cap U\neq \varnothing\}}  |\psi(g_U)f(g_U)\xi(g_U^{-1}g)|^2 \\ 
&\quad\leq \sup_{g\in G}|\psi(g)|^2\sup\Big\{|f(g)|^2\sum_{g\in s^{-1}(x)}|\xi(g)|^2~\Big|~g\in \bigcup_{s(h)=x} s^{-1}(r(h))\Big\} \\
& \quad= \sup_{g\in G}|\psi(g)|^2\sup\Big\{|f(g)|^2\|\xi\|^2~\Big|~g\in \bigcup_{s(h)=x} s^{-1}(r(h))\Big\}.
\end{align*}
As the expression $\sup\{|f(g)|^2~|~g\in \cup_{s(h)=x} s^{-1}(r(h))\}$ is easily seen to be a lower bound for $\|\pi_x(f)\|$, we are done. 
\end{proof}

\begin{lemma}\label{commlem}
For any $\epsilon>0$ and compact subset $K$ of $G$, there exists $\delta>0$ such that if $\phi\in C_c(G^{(0)})$ satisfies 
$$
\sup_{g\in K} |\phi(r(g))-\phi(s(g))|<\delta,
$$
then the commutator $[f,\phi]$ has norm at most $\epsilon\|\phi\|\|f\|$ for any $f\in C_K(G)$.
\end{lemma}

\begin{proof}
With assumptions as above, the commutator is the element of $C_c(G)$ given by
$$
[f,\phi](g)=f(g)\phi(s(g))-\phi(r(g))f(g)
$$
in other words, it is the function $(\phi\circ s-\phi\circ r)\cdot f$, where `$\cdot$' denotes pointwise multiplication.  Fix now a compact set $K'$ containing an open neighbourhood of $K$, and let $\psi$ be any function that agrees with $\phi\circ s-\phi\circ r$ on $K$, vanishes outside $K'$ and is bounded above by $\sup_{g\in K}|(\phi\circ s-\phi\circ r)(g)|$.  Lemma \ref{multlem} then implies that there exists $M$ such that for any $f\in C^*_K(G)$
\begin{align*}
\|[f,\phi]\|_{C^*_r(G)}&=\|\psi\cdot f\|_{C^*_r(G)}\leq M\sup_{g\in G}|\psi(g)|\|f\|_{C^*_r(G)} \\ &\leq M\sup_{g\in K}|(\phi\circ s-\phi\circ r)(g)|\|f\|_{C^*_r(G)};
\end{align*}
taking $\delta=\epsilon/M$ thus works.
\end{proof}

\begin{proof}[Proof of Theorem \ref{fnd}]
Let a finite subset $\mathcal{F}$ of $C^*_r(G)$ and $\epsilon>0$ be given.  As $C_c(G)$ is dense in $C^*_r(G)$ we may assume that $\mathcal{F}$ is a subset of $C_c(G)$; let $K$ be a compact subset of $G$ such that $K=K^{-1}$ and such that each element $f\in \mathcal{F}$ is supported in $K$.  Let $\widetilde{C^*_r(G)}$ denote the unitization of $C^*_r(G)$, and let $\widetilde{C^*_K(G)}$ denote the closed subspace of $\widetilde{C^*_r(G)}$ spanned by $C^*_K(G)$ and the identity of $\widetilde{C^*_r(G)}$; as $K=K^{-1}$, this is an operator subsystem of $\widetilde{C^*_r(G)}$.  Using Proposition \ref{dadpou} (see also Remark \ref{str rem} part \eqref{dadpou still good}), there exists a collection $\{\phi_i:G^{(0)}\to [0,1]~|~i\in\{0,...,d\}\}$ of continuous compactly supported functions on $G^{(0)}$ with the following properties:
\begin{enumerate}[(i)]
\item for each $x\in r(K)\cup s(K)$ we have $\sum_{i=0}^d \phi_i^2(x)=1$;
\item for each $i$ there is a small open subgroupoid $H_i$ of $G$ such that for any $f\in \widetilde{C^*_K(G)}$, the element $\phi_if\phi_i$ is in contained in the sub-$C^*$-algebra $C^*_r(H_i)$ of $C^*_r(G)$;
\item for each $i$ and $g\in K$,
$$
|\phi_i(r(g))-\phi_i(s(g))|<\frac{\epsilon}{2(d+1)}.
$$ 
\end{enumerate}
From the third point above and Lemma \ref{commlem} we have
\begin{equation}\label{secondinq}
\|[f,\phi_i]\|<\frac{\epsilon}{2(d+1)}\|f\|
\end{equation}
for any $f\in C^*_K(G)$, so in particular for all $f\in \mathcal{F}$.   

Now, for each $i\in\{0,...,d\}$, the formula
$$
\Phi_i:\widetilde{C^*_K(G)}\to C^*_r(H_i),~~~f\mapsto \phi_i f\phi_i 
$$
defines a contractive completely positive map.   Let 
$$
\Psi_i:C^*_r(H_i)\to C^*_r(G)
$$
be the canonical inclusion $*$-homomorphism.  For each $i\in \{0,...d\}$ we now have a triangle of maps
\begin{equation}\label{first tri}
\xymatrix{ \widetilde{C^*_{K}(G)} \ar[dr]^{\Phi_i} & & C^*_{r}(G) \\ & C^*_r(H_i) \ar[ur]^{\Psi_i} & },
\end{equation}
where each $\Phi_i$ is contractive and completely positive, and each $\Psi_i$ is a $*$-homomorphism, so in particular order zero.  Moreover, for every $f\in \mathcal{F}$  
\begin{align*}
\sum_{i=0}^d\Psi_i(\Phi_i(f)) &=\sum_{i=0}^d\phi_if\phi_i =\sum_{i=0}^d\phi_i^2f+\sum_{i=0}^d\phi_i[f,\phi_i]=f+\sum_{i=0}^d\phi_i[f,\phi_i],
\end{align*}
whence line \eqref{secondinq} implies that
$$
\Big\|\sum_{i=0}^d\Psi_i(\Phi_i(f))-f\Big\|<\epsilon/2
$$
for all $f\in \mathcal{F}$.

On the other hand, we have already shown that each $C^*_r(H_i)$ has nuclear dimension at most the covering dimension of $H_i^{(0)}$ in Proposition \ref{smallprop}, and this is at most $N$ as covering dimension does not increase under taking open subsets (this follows for example from \cite[Proposition 2.5]{Winter:2010eb}).  Combining the triangle in line \eqref{first tri} with approximations to each $C^*_r(H_i)$ arising from the definition of nuclear dimension gives triangles
$$
\xymatrix{ \widetilde{C^*_{K}(G)} \ar[dr]^{\Phi_{ij}} & & C^*_{r}(G) \\ & F_{ij} \ar[ur]^{\Psi_{ij}} & }, \quad i\in \{0,...d\},~j\in \{0,...,N\}
$$
where each $\Phi_{ij}$ is contractive and completely positive, each $\Psi_{ij}$ is contractive and order zero, each $F_{ij}$ is a finite dimensional $C^*$-algebra, and 
$$
\Big\|\sum_{i=0}^d \sum_{j=0}^N\Psi_{ij}(\Phi_{ij}(f))-f\Big\|<\epsilon
$$
for all $f\in \mathcal{F}$.  Finally, the finite dimensional version of Arveson's extension theorem (see for example \cite[Section 1.6]{Brown:2008qy}) implies that each $\Phi_{ij}$ extends to a contractive completely positive map
$$
\Phi_{ij}:\widetilde{C^*_r(G)}\to F_{ij},
$$
and restricting each of these maps to $C^*_r(G)$ gives approximating triangles 
$$
\xymatrix{ C^*_r(G) \ar[dr]^{\Phi_{ij}} & & C^*_{r}(G) \\ & F_{ij} \ar[ur]^{\Psi_{ij}} & }, \quad i\in \{0,...d\},~j\in \{0,...,N\}
$$
as required by the definition of nuclear dimension $(d+1)(N+1)-1$.
\end{proof}

\subsection*{Consequences}

Finally, we spell out a few consequences.  Many of these are known results, but we think that combining the proofs under one common dynamical framework has some interest.

The first result is much easier to check directly! See Winter \cite[Remark 2.2 (iii)]{Winter:2010eb} and Kirchberg-Winter \cite[Example 4.1]{Kirchberg:2004uq}; nonetheless, it seemed interesting that it fits directly into our framework.

\begin{corollary}\label{affnd}
Separable AF $C^*$-algebras have nuclear dimension zero.
\end{corollary}

\begin{proof}
Renault \cite[Proposition III.1.15]{Renault:1980fk} has shown that any (separable) AF $C^*$-algebra arises as the $C^*$-algebra of a locally finite groupoid (in the sense of Example \ref{lfgpd}) with unit space of covering dimension zero.  The result thus follows directly from Example \ref{lfgpd} and Theorem \ref{fnd}.
\end{proof}

The next corollary is due to Winter-Zacharias: see  \cite[Section 8]{Winter:2010eb}.

\begin{corollary}\label{roefnd}
Let $X$ be a bounded geometry coarse space of asymptotic dimension $d$.  Then the uniform Roe algebra $C^*_u(X)$ has nuclear dimension at most $d$. 
\end{corollary}

\begin{proof}
The uniform Roe algebra of $X$ is naturally isomorphic to the reduced groupoid $C^*$-algebra of the coarse groupoid $G(X)$: see \cite[Proposition 10.29]{Roe:2003rw} for a proof.  The result now follows on combining Theorem \ref{fadthe}, Theorem \ref{fnd} (see also Remark \ref{str rem}, part \eqref{fad str}), and the fact that the unit space of $G(X)$ is $\beta X$, which has covering dimension zero.
\end{proof}

The next corollary is due to Toms-Winter \cite[Section 3]{Toms:2013vn}.

\begin{corollary}\label{sysfnd}
Let $\Z\lefttorightarrow X$ be a minimal action of $\Z$ on an infinite second countable compact space $X$ of covering dimension $N$.  Then the nuclear dimension of $C(X)\rtimes_r\Z$ is at most $2N+1$.
\end{corollary}

\begin{proof}
Combine Theorem \ref{minz}, Lemma \ref{gpgpd}, Theorem \ref{fnd}, and the well-known (and easily checked) fact that the reduced groupoid $C^*$-algebra of a transformation groupoid $\Gamma\ltimes Y$ identifies naturally with $C(Y)\rtimes_r\Gamma$.
\end{proof}

It is worth noting explicitly that the proofs of Corollary \ref{roefnd} and Corollary \ref{sysfnd} in the original references given above are quite similar to our proof of Theorem \ref{fnd}: all involve constructing `almost invariant' partitions of unity, and using these to `cut down' to subhomogeneous $C^*$-algebras whose nuclear dimension can be directly estimated.  This only became apparent to us `after the fact', but we hope it helps to clarify how this style of argument that estimates nuclear dimension from geometric and / or dynamic assumptions is built up.

The following result seems to be new as stated.  Note however that if $\Gamma\lefttorightarrow X$ is a minimal, free, amenable action of a \emph{non-amenable} group on the Cantor set, then $C(X)\rtimes_r\Gamma$ is a Kirchberg algebra in the UCT class; the main result of \cite{Ruiz:2014aa} (see also \cite{Bosa:2014zr}) thus implies that the nuclear dimension of $C(X)\rtimes_r\Gamma$ is one.  As any exact (non-amenable) group admits such an action (see \cite[Theorem 6.11]{Rordam:2010kx}), and as finite asymptotic dimension implies exactness (see \cite[Lemma 4.3]{Higson:2000dp} and \cite{Ozawa:2000th}) but not conversely, the corollary below is not optimal, neither with respect to the class of groups covered, nor with respect to the estimate on nuclear dimension.

\begin{corollary}\label{cantfnd}
Let $\Gamma$ be a countable discrete group with finite asymptotic dimension.  Then $\Gamma$ admits a minimal, free action on the Cantor set $X$ such that $C(X)\rtimes_r\Gamma$ has nuclear dimension bounded above by the asymptotic dimension of $\Gamma$.
\end{corollary}

\begin{proof}
Combine Lemma \ref{gpgpd}, Theorem \ref{mincor}, Theorem \ref{fnd}, and again that the reduced groupoid $C^*$-algebra of a transformation groupoid $\Gamma\ltimes Y$ identifies naturally with $C(Y)\rtimes_r\Gamma$.
\end{proof}

We should note that there are many other results in the literature that give estimates on nuclear dimension or other $C^*$-algebraic regularity properties based on conditions on actions: see for example \cite{Szabo:2015mw,Hirshberg:2015hb,Szabo:2014aa,Hirshberg:2015aa,Elliott:2014aa}.  Many of these results go further than ours in at least some ways: for example, \cite{Szabo:2015mw} proves fairly general finite dimensionality results for $\Z^n$-actions, \cite{Hirshberg:2015aa} treats non-free $\Z$-actions, \cite{Elliott:2014aa} treats some $\Z$-actions on non-finite dimensional spaces, and several works deal with some actions on noncommutative $C^*$-algebras.  It would be interesting to clarify the relationships holding between the various conditions involved in these results and ours.

The final result is not strictly a `corollary' as such, but follows from exactly the same method of proof; it does not require any separability assumptions as these were only used in the above in the appeal to Theorem \ref{shdr}, and to avoid complications from general topology with respect to dimension theory.  It is no doubt possible to prove it more directly.

\begin{corollary}\label{dadamen}
Let $G$ be an \'{e}tale groupoid with compact unit space, and with finite dynamic asymptotic dimension.  Then $G$ is amenable. 
\end{corollary}

\begin{proof}
The method of proof of Theorem \ref{fnd} implies in particular that $C^*_r(G)$ is nuclear (whether or not the unit space of $G$ is finite dimensional).  Hence $G$ is amenable by \cite[Theorem 5.6.18]{Brown:2008qy}.
\end{proof}

\bibliography{Generalbib}

\end{document}